\documentclass[8pt,nonatbib]{elsarticle}
\usepackage[utf8]{inputenc}
\usepackage{amsmath,amssymb,amsthm}
\usepackage{latexsym}            
\usepackage{amsthm}
\usepackage{booktabs}
\newtheorem{remark}{Remark}
\usepackage{epsfig}
\usepackage{epstopdf}
\usepackage{amsfonts,amscd, makecell}
\usepackage[ruled,vlined]{algorithm2e}
\usepackage{color,comment}
\usepackage{amsbsy,bm}
\usepackage{caption}
\usepackage{subcaption}
\usepackage{setspace}
\usepackage{enumitem}
\usepackage{xcolor}
\usepackage[normalem]{ulem}
\usepackage{placeins}

\newlist{steps}{enumerate}{1}
\setlist[steps]{label=Step \arabic*:, 
                leftmargin=2.5cm,   
                rightmargin=0.8cm,
                labelwidth=1.5cm,     
                align=left,         
                itemindent=0pt,     
                before={\setlength{\parindent}{0pt}}} 

\usepackage[utf8]{inputenc}          
\usepackage[T1]{fontenc}

\newcommand\vel{\mathbf{u}}

\newcommand\bv{\mathbf{v}}
\newcommand\bx{\boldsymbol{x}}

\newcommand\bn{\mathbf{n}}
\newcommand\bX{\mathbf{X}}

\newcommand\bF{\mathbf{F}}

\newcommand\bD{\boldsymbol{D}}

\newcommand{\blue}[1]{#1}

\newcommand{\be}{\begin{eqnarray}}
\newcommand{\ee}{\end{eqnarray}}
\newcommand{\ben}{\begin{eqnarray*}}
\newcommand{\een}{\end{eqnarray*}}

\newcommand {\na} {{\nabla}}

 \allowdisplaybreaks \allowdisplaybreaks[4]
\makeatletter
\let\c@author\relax
\makeatother
\usepackage[backend=biber,style=numeric,sorting=nyt]{biblatex}
\usepackage[colorlinks=true,linkcolor=blue,citecolor=blue,urlcolor=blue]{hyperref}
\begin{document}
\begin{frontmatter}
  \title{A quasi-monolithic localized high-order ALE finite element method for multi-scale fluid-structure interaction problems}

  \author[1]{Lingyue Shen}
  \ead{lingyue.shen@raymind.com}
  \author[2,3]{Qi Xin}
  \ead{qixin1@link.cuhk.edu.cn}
  \author[1]{Yan Chen}
  \ead{yanc@raymind.com}
  \author[1]{Jiarui Han \corref{cor1}}
  \ead{han@raymind.com}
  \author[1]{Yumiao Zhang}
  \ead{yumiao.zhang@raymind.com}
  \author[5]{Jinchao Xu}
  \ead{jinchao.xu@kaust.edu.sa}
  \author[2,3,4]{Shihua Gong \corref{cor1}}
  \ead{gongshihua@cuhk.edu.cn}

  \affiliation[1]{organization={Shenzhen Raymind Biotechnology Co., Ltd.},
    addressline={Shenzhen},
    postcode={518129},
    state={Guangdong},
    country={China}}

  \affiliation[2]{organization={School of Science and Engineering, The Chinese University of Hong Kong},
    addressline={Shenzhen},
    postcode={518000},
    state={Guangdong},
    country={China}}

  \affiliation[3]{organization={Shenzhen International Center for Industrial and Applied Mathematics, Shenzhen Research Institute of Big Data},
    addressline={Shenzhen},
    postcode={518000},
    state={Guangdong},
    country={China}}
    \affiliation[4]{
    organization={Shenzhen Loop Area  Institute},
    addressline={Shenzhen},
    postcode={518000},
    state={Guangdong},
    country={China}}

  \affiliation[5]{organization={Applied Mathematics and Computational Sciences, CEMSE Division, King Abdullah University of Science and Technology},
    city={Thuwal},
    postcode={23955},
    country={Saudi Arabia}}

  \cortext[cor1]{Corresponding authors}

  \date{}


  \begin{abstract}

    This paper presents a quasi-monolithic localized high-order arbitrary Lagrangian-Eulerian (qMLH-ALE) finite element method for multi-scale fluid-structure interaction (FSI) in microfluidic systems. The fluid momentum, the incompressible Neo-Hookean constitutive law, and the left Cauchy-Green tensor $\mathcal{B}$ are assembled into a single implicit system, while the harmonic mesh extension is updated explicitly in a staggered manner. Isoparametric $\mathcal{P}_2$ elements provide third-order geometric approximation of curved fluid-solid interfaces, and a second-order implicit-explicit partitioned Runge-Kutta scheme delivers second-order temporal accuracy without the dissipation of backward Euler. A localized updating strategy confines the moving mesh and the deformation history to a body-fitted sub-domain coupled with a precomputed steady background flow, bridging the scale disparity between local FSI dynamics and the macroscopic microchannel geometry. The Turek-Hron FSI3 benchmark, performed at unit fluid-solid density ratio, reproduces the reference beam-tip amplitude and frequency within $3\%$, confirming stability under the strong added-mass coupling that destabilizes conventional partitioned schemes. Three-dimensional particle-focusing simulations in spiral microchannels further illustrate the framework on long-range multi-scale problems.
  \end{abstract}

\end{frontmatter}

\section{Introduction}
Fluid-structure interaction (FSI) computations underpin engineering design across aerospace, naval architecture, and biomedicine: aerodynamic loading on aircraft wings and rotor blades \cite{hu2010material, joshi2020variational, gupta2024modeling, jayakumar2024design}, hull-ice and hull-wave impact in shipbuilding \cite{song2017fluid, lee2017full, rath2025phase}, and blood-cell dynamics in vascular flow \cite{almomani2012sharp, marth2016margination, xiao2017effects, syed2023modeling}.

In recent years, research on biological microfluidic devices has grown rapidly \cite{ozbey2016inertial, erdem2020differential, dalili2019review, nasiri2020microfluidic}. These devices control particle trajectories and dynamical states by manipulating the microscale flow field within the chip, with applications in medical diagnostics such as cancer-cell separation \cite{shen2024spiral, horade2023analysis, betyar2026numerical} and single-cell encapsulation \cite{gardner2022deep, zhao2021injectable}. The performance of these devices depends on the accuracy of the underlying flow and particle dynamics, and FSI simulation is therefore directly tied to design efficiency and performance optimization. Three structural obstacles, however, restrict the applicability of existing FSI solvers to microfluidic problems. The first is the multiscale gap between local FSI phenomena (tens of micrometers, low Reynolds number) and the macroscopic device dimensions (several centimeters); a global FSI computation over the entire domain incurs prohibitive overhead, and periodic boundary conditions are excluded in asymmetric devices such as spiral channels \cite{shen2024spiral}. The second is the long-range migration of particles, which requires sustained high-fidelity tracking of the FSI dynamics over distances orders of magnitude larger than the particle itself. The third is the geometric representation of the fluid-structure interface, since devices with microscale sub-structures such as deterministic lateral displacement (DLD) arrays \cite{liu2021cascaded, jiang2025chip} are sensitive to small variations in particle and obstacle curvature. Existing FSI solvers address these obstacles in isolation rather than jointly.


Computing methods developed for general FSI problems fall into three families, each handling the three obstacles above with a different trade-off. Immersed boundary methods (IBM) on Cartesian grids \cite{peskin1977numerical, sotiropoulos2014immersed} keep the mesh fixed during simulation, giving high efficiency and good scalability on structured-grid parallel architectures, but their stair-cased representation of curved boundaries compromises geometric fidelity for microscale sub-structures \cite{powar2025recent}.

Diffuse interface methods such as phase-field and level-set formulations \cite{liu2001eulerian, hong2021hybrid, rath2021interface, mao20243d} replace the discrete force coupling of IBM with a smooth order parameter, enforcing thermodynamic consistency and enabling a unified monolithic discretization of the fluid-solid continuum. The smooth transition, however, smears the physical interface over several grid cells, again limiting the accuracy of wall-shear and contact predictions on microscale geometries.

  {The Arbitrary Lagrangian-Eulerian (ALE) method \cite{hirt1974arbitrary} and the Deforming-Spatial-Domain/Stabilized Space-Time (DSD/SST) method \cite{tian2015fsi} take the opposite stance: a dynamic body-fitted mesh moves with the structural motion, so the computational nodes align exactly with the physical boundary \cite{tian2015fsi, dehghani2023finite, li2025numerical}. Interfacial stresses and wall-shear forces are recovered with high resolution, addressing the geometric-fidelity obstacle directly. The remaining bottleneck is mesh management: high-fidelity tracking relies on robust mesh-update or local-remeshing strategies to avoid element inversion under large structural deformation and long-range migration, which is precisely the regime imposed by the first two obstacles in microfluidic problems.}

  {Based on the aforementioned research and methodology, this study introduces a quasi-monolithic localized sharp-interface FSI framework built on the ALE description, designed around three targeted contributions. First, the fluid momentum, the incompressible Neo-Hookean constitutive equation, and the evolution of the left Cauchy-Green tensor $\mathcal{B}$ are assembled into a single implicit system, while the harmonic mesh extension is updated explicitly in a staggered manner; this fluid-solid quasi-monolithic formulation maintains strong interface coupling under matched fluid-solid density, where the added-mass effect destabilizes conventional partitioned schemes. Second, isoparametric $\mathcal{P}_2$ elements on body-fitted curved meshes provide third-order $\mathcal{O}(h^3)$ geometric approximation of the fluid-solid interface, removing the staircase artifacts of fixed Cartesian or immersed-boundary discretizations and recovering the wall-shear distribution that governs separation efficiency in DLD and spiral microfluidic devices. Third, a localized updating strategy exploits the exponentially decaying hydrodynamic perturbation in low-Reynolds-number confined flow: the moving mesh and the structural deformation history are restricted to a body-fitted sub-domain that follows the particle, while the surrounding flow is read from a precomputed steady background field. This decomposition replaces a global remeshing of the centimeter-scale channel with a sub-domain solve, making three-dimensional long-range multi-scale FSI tractable on standard workstations.}

  Two additional choices realize these contributions in a numerically stable manner. To attain second-order temporal accuracy without the dissipation of backward Euler, the ALE-based weak form is pulled back to a fixed reference configuration so that the moving-domain problem reduces to a stationary-domain ODE system, on which an implicit-explicit partitioned Runge-Kutta (IMEX-PRK) scheme \cite{kennedy2003additive} preserves the second-order rate while handling the disparate stiffness scales of FSI. Spatial discretization is performed by the finite element method on unstructured meshes with isoparametric $\mathcal{P}_2$ elements, which deliver the $\mathcal{O}(h^3)$ geometric fidelity claimed above.

The paper is arranged as follows: Section 2 derives the FSI model and its variational formulation based on the ALE framework. Section 3 presents the high-order spatial and temporal discretization of the derived problem. Section 4 states a local updating algorithm suitable for large {displacement} and multiscale FSI simulation. Section 5 validates the model through benchmark tests and provides 2D/3D numerical demonstrations. {Section 6 concludes the paper.}

\section{Arbitrary Lagrangian Eulerian (ALE) Formulation of a Sharp Interface FSI Model}

In this section, we first introduce an incompressible elastic model in fluid structure interaction problems, then derive a sharp interface formulation and the corresponding governing equations.

\subsection{Incompressible elastic model for fluid-structure interaction}
The FSI system under investigation couples an incompressible Newtonian fluid with an elastic structure governed by the incompressible Neo-Hookean hyperelastic model. The entire computational domain is denoted by $\Omega \subset \mathbb{R}^d$ ($d=2,3$), which is a {fixed, bounded region}. At time $t$, the solid occupies the subdomain $\Omega_s^t \subset \Omega$, which is the image of the reference configuration $\hat{\Omega}_s$ under the deformation map $\boldsymbol{X}(\cdot,t): \hat{\Omega}_s \rightarrow \Omega_s^t$. The fluid occupies the remaining part of the domain, denoted as $\Omega_f^t = \Omega \setminus \overline{\Omega_s^t}$. The interface between the fluid and solid is represented by $\Gamma^t = \partial \Omega_s^t \cap \partial \Omega_f^t$. {We assume that the solid is always fully immersed in the fluid and does not contact the external boundary}, i.e., $\partial \Omega_s^t \cap \partial \Omega = \emptyset$ for all $t \in [0,T]$. Consequently, the fixed boundary $\partial \Omega$ is in contact with the fluid only.

The motion of the fluid is governed by the incompressible Navier-Stokes equations:
\be
\label{eq:NS}\left\{\begin{array}{ll}
  \frac{\partial \vel_f}{\partial t}+(\vel_f\cdot \na)\vel_f=\nabla\cdot \boldsymbol{\sigma}_f~, & \mbox{in~}\Omega_f^t~, \\[3mm]
  \na\cdot \vel_f=0~,                                                                            & \mbox{in~}\Omega_f^t~,
\end{array}\right.
\ee
where $\vel_f$\ denote the fluid velocity. The fluid stress tensor $\boldsymbol{\sigma}_f$ is given by
\be
\boldsymbol{\sigma}_f = -P_f \mathbf{I} + \frac2{Re} \bD(\vel_f)~,
\ee
with $\bD(\vel_f) = \frac12(\na \vel_f + (\na \vel_f)^T)$ representing the rate-of-deformation tensor, $P_f$ the fluid pressure, and $Re$ being the Reynolds number of fluid.

For the structure, we adopt a neo-Hookean elastic model. The motion is described by the following equations:
\be
\label{eq:solid}\left\{\begin{array}{ll}
  \partial_t \vel_s + (\vel_s\cdot\nabla )\vel_s = \nabla \cdot \boldsymbol{\sigma}_s~, & \mbox{in~}\Omega_s^t~, \\[3mm]
  \nabla\cdot\vel_s = 0~,                                                               & \mbox{in~}\Omega_s^t~,
\end{array}\right.
\ee
where $\vel_s = \partial_t \bX \circ \bX^{-1}$ is the solid velocity. The solid stress tensor $\boldsymbol{\sigma}_s$ is defined as
\be
\boldsymbol{\sigma}_s = -P_s \mathbf{I} + \mu_s(\mathcal{B}-\mathbf{I})~,
\ee
where $P_s$ is the solid pressure, and $\mu_s$ is the dimensionless shear modulus, defined as $\mu_s^* = \mu_s / (\rho_f U^2)$ under the inertial scaling adopted here, so that the momentum balance is non-dimensionalized by the dynamic pressure $\rho_f U^2$ rather than by $Re$. We retain the symbol $\mu_s$ for $\mu_s^*$ in the remainder of the paper. The left Cauchy-Green deformation tensor is $\mathcal{B} = \bF \bF^T$, where $\bF = \na_{\bX} \bX$ is the deformation gradient. For the strictly incompressible case ($\nu_s = 1/2$), the dimensional shear modulus relates to the Young's modulus through $\mu_s^{\mathrm{dim}} = E / 3$. $\bF$ satisfies the transport equation \cite{liu2001eulerian}
\[
  \frac{\partial \bF}{\partial t} + (\vel_s \cdot \na) \bF = (\na \vel_s)^T \bF,
\]
leading to the evolution equation for $\mathcal{B}$
\[
  \frac{\partial \mathcal{B}}{\partial t} + (\vel_s \cdot \na) \mathcal{B} - (\na \vel_s)^T \mathcal{B} - \mathcal{B} (\na \vel_s) = 0,
\]
which expresses the vanishing of the upper-convected time derivative of $\mathcal{B}$.

The fluid and structure are coupled through the following interface conditions on $\Gamma^t$:
\be
\label{eq:interface-conditions}\left\{\begin{array}{ll}
  \vel_f = \vel_s~,                                        & \mbox{on~}\Gamma^t~, \\[3mm]
  \boldsymbol{\sigma}_f  \bn = \boldsymbol{\sigma}_s \bn~, & \mbox{on~}\Gamma^t~,
\end{array}\right.
\ee
where $\bn$ is the unit normal vector on $\Gamma^t$ pointing from the solid to the fluid. Since the velocities are continuous across the interface . Thus, together with suitable boundary and initial conditions, the coupled FSI system is given by:
\be
\label{eq:FSI}\left\{\begin{array}{ll}
  \frac{\partial \vel_s}{\partial t}+(\vel_s\cdot \na)\vel_s+\na P_s=\mu_s\na \cdot (\mathcal{B}-\mathbf{I})~, & \mbox{in~}\Omega_s^t~,     \\[3mm]
  \frac{\partial \vel_f}{\partial t}+(\vel_f\cdot \na)\vel_f+\na P_f=\frac{1}{Re}\na\cdot(2\bD(\vel_f))~,                         & \mbox{in~}\Omega_f^t~,     \\[3mm]
  \na\cdot \vel_s=0~,                                                                                                             & \mbox{in~}\Omega_s~,       \\[3mm]
  \na\cdot \vel_f=0~,                                                                                                             & \mbox{in~}\Omega_f~,       \\[3mm]
  \frac{\partial\mathcal{B}}{\partial t}+(\vel_s\cdot\na)\mathcal{B} - (\na \vel_s)^T \mathcal{B} - \mathcal{B} \na \vel_s = 0~,  & \mbox{in~}\Omega_s^t~,     \\[3mm]
  \vel_f = \vel_s~,                                                                                                               & \mbox{on~}\Gamma^t~,       \\[3mm]
  \boldsymbol{\sigma}_f \bn = \boldsymbol{\sigma}_s \bn~,                                                                         & \mbox{on~}\Gamma^t~,       \\[3mm]
  \vel_s = \vel_f = 0~,                                                                                                           & \mbox{on~}\partial\Omega~, \\[3mm]
  \vel_{s}(\bx,0) = \vel_{s, 0}(\bx)~,\vel_{f}(\bx,0) = \vel_{f, 0}(\bx)~,~\mathcal{B}(\bx,0) = \mathbf{I}~,                      & \mbox{in~}\Omega~,
\end{array}\right.
\ee

\begin{remark}
  The left Cauchy-Green tensor $\mathcal{B}$ is not an independent unknown; it is the auxiliary kinematic field constrained by $\mathcal{B} = \mathbf{F}\mathbf{F}^\top$. We retain $\mathcal{B}$ as a primary variable and advance it through its own evolution equation for three reasons. First, separating $\mathbf{F}\mathbf{F}^\top$ from the stress modularizes the Newton linearization of the nonlinear Neo-Hookean term. Second, evolving $\mathcal{B}$ directly avoids reconstructing the deformation gradient from the displacement, which would amplify numerical noise on ALE meshes undergoing severe distortion. Third, treating $\mathcal{B}$ as a state field on the Eulerian-ALE configuration enables consistent high-order projection of the deformation history during the local-mesh transfers introduced in Section 4.
\end{remark}

\subsection{Variational formulation and ALE mapping}
Let us denote by $H^1(\Omega)$ the Sobolev space of  functions whose first derivatives and themselves are square-integrable, and by $L^2(\Omega)$ the space of square-integrable functions.

To derive a natural weak formulation, we first introduce the following function spaces for the separated fluid and structural domains respectively:

\begin{align*}
  \tilde{\mathbf{V}}^t & := \{(\mathbf{v}_f, \mathbf{v}_s)~ |~ \mathbf{v}_f \in H^1(\Omega_f^t)^d, \mathbf{v}_s \in H^1(\Omega_s^t)^d; \mathbf{v}_f = \mathbf{v}_s \text{ on } \Gamma^t, \mathbf{v}_f = 0 \text{ on } \partial\Omega_f^t \setminus \Gamma^t, \mathbf{v}_s = 0 \text{ on } \partial\Omega_s^t \setminus \Gamma^t \}~, \\
  \tilde{W}^t          & := \{(q^f, q^s) ~| ~q^f \in L^2(\Omega_f^t), q^s \in L^2(\Omega_s^t); \int_{\Omega_f^t} q^f d\bx = 0, \int_{\Omega_s^t} q^s d\bx = 0 \}~,                                                                                                                                                                   \\
  \mathbb{B}^t         & := H^1(\Omega_s^t)^{d \times d}~.\end{align*}

Multiplying the fluid and solid momentum equations in \eqref{eq:FSI} by test functions $\tilde{\mathbf{v}}_f$ and $\tilde{\mathbf{v}}_s$ respectively, integrating by parts over $\Omega_f^t$ and $\Omega_s^t$, and summing them up, we obtain the separated weak formulation: find $[(\mathbf{u}_f, \mathbf{u}_s), (P^f, P^s), \mathcal{B}] \in \tilde{\mathbf{V}}^t \times \tilde{W}^t \times \mathbb{B}^t$ such that for any test functions $[(\tilde{\mathbf{v}}_f, \tilde{\mathbf{v}}_s), (\tilde{q}^f, \tilde{q}^s), \mathcal{W}] \in \tilde{\mathbf{V}}^t \times \tilde{W}^t \times \mathbb{B}^t$,

\begin{equation}
  \begin{aligned}
    \label{eq:separated_weak}
     & \left(\frac{\partial \mathbf{u}_f}{\partial t}, \tilde{\mathbf{v}}_f \right)_{\Omega_f^t} + \left((\mathbf{u}_f \cdot \nabla)\mathbf{u}_f, \tilde{\mathbf{v}}_f \right)_{\Omega_f^t} + \frac{1}{Re} \left(2\mathbf{D}(\mathbf{u}_f), \mathbf{D}(\tilde{\mathbf{v}}_f) \right)_{\Omega_f^t} - \left(P^f, \nabla \cdot \tilde{\mathbf{v}}_f \right)_{\Omega_f^t} \\
     & + \left(\frac{\partial \mathbf{u}_s}{\partial t}, \tilde{\mathbf{v}}_s \right)_{\Omega_s^t} + \left((\mathbf{u}_s \cdot \nabla)\mathbf{u}_s, \tilde{\mathbf{v}}_s \right)_{\Omega_s^t} - \left(P^s, \nabla \cdot \tilde{\mathbf{v}}_s \right)_{\Omega_s^t} + \mu_s \left(\mathcal{B}-\mathbf{I}, \nabla \tilde{\mathbf{v}}_s \right)_{\Omega_s^t}       \\
     & = \int_{\Gamma^t} (\boldsymbol{\sigma}_f \mathbf{n}_f) \cdot \tilde{\mathbf{v}}_f dS + \int_{\Gamma^t} (\boldsymbol{\sigma}_s \mathbf{n}_s) \cdot \tilde{\mathbf{v}}_s dS = 0~,                                                                                                                                                                                \\
     & \left(\nabla \cdot \mathbf{u}_f, \tilde{q}^f \right)_{\Omega_f^t} = 0, \quad \left(\nabla \cdot \mathbf{u}_s, \tilde{q}^s \right)_{\Omega_s^t} = 0~,                                                                                                                                                                                                           \\
     & \left(\frac{\partial \mathcal{B}}{\partial t}, \mathcal{W} \right)_{\Omega_s^t} + \left((\mathbf{u}_s \cdot \nabla)\mathcal{B}, \mathcal{W} \right)_{\Omega_s^t} - \left((\nabla \mathbf{u}_s)^T \mathcal{B}, \mathcal{W} \right)_{\Omega_s^t} - \left(\mathcal{B} \nabla \mathbf{u}_s, \mathcal{W} \right)_{\Omega_s^t} = 0~.
  \end{aligned}
\end{equation}

Note that the interface integral terms on the right-hand side of \eqref{eq:separated_weak} naturally vanish due to the traction balance condition $\boldsymbol{\sigma}_f \mathbf{n}_f + \boldsymbol{\sigma}_s \mathbf{n}_s = 0$ and the kinematic continuity $\tilde{\mathbf{v}}_f = \tilde{\mathbf{v}}_s$ on $\Gamma^t$ built into the space $\tilde{\mathbf{V}}^t$.

Because $\Omega = \Omega_f^t \cup \Omega_s^t$, we can define the unified velocity $\mathbf{u}$ over the entire domain $\Omega$ by smoothly extending the piece-wise definitions. Consequently, we define the following global function spaces:

\begin{align*}
  \mathbf{V}^t & =\{\mathbf{v}\in H^1(\Omega)^d|~\mathbf{v}=0~\mbox{on}~\partial\Omega\}~, \\
  Q_s^t        & =\{q^s\in L^2(\Omega_s^t)| \int_{\Omega_s^t}q^s d\bx = 0\}~,              \\
  Q_f^t        & =\{q^f\in L^2(\Omega_f^t)| \int_{\Omega_f^t}q^f d\bx = 0\}~.              \\
\end{align*}


Utilizing these global variables and the inner product notations $(\cdot,\cdot)_{\Omega}$, $(\cdot,\cdot)_{\Omega_s^t}$, and $(\cdot,\cdot)_{\Omega_f^t}$, we are able to rewrite \eqref{eq:separated_weak} as the following equivalent monolithic weak formulation: find $(\mathbf{u},P^{s},P^{f},\mathcal{B})\in \mathbf{V}^t\times Q^t_s\times Q^t_f\times\mathbb{B}^t$ such that for all $(\bv,q^s,q^f,\mathcal{W})\in \mathbf{V}^t\times Q^t_s\times Q^t_f\times\mathbb{B}^t$

\be\left\{\begin{array}{ll}
  \label{weakform}
  (\frac{\partial \vel}{\partial t}, \bv)_\Omega +((\vel\cdot \na)\vel,\bv)_\Omega + \frac{1}{Re}(2\mathbf{D}(\mathbf{u}), \mathbf{D}(\bv))_{\Omega_f} \\
  ~~~~~~~~-(P^{s}, \na\cdot \bv)_{\Omega_s^t} -(P^{f}, \na\cdot \bv)_{\Omega_f^t} +\mu_s((\mathcal{B}-\mathbf{I}), \na\bv)_{\Omega_s^t}=0,      \\[3mm]
  (\na\cdot \vel,q^s)_{\Omega_s^t}=0~,                                                                                                                 \\[3mm]
  (\na\cdot \vel,q^f)_{\Omega_f^t}=0~,                                                                                                                 \\[3mm]
  (\frac{\partial \mathcal{B}}{\partial t},\mathcal{W})_{\Omega_s^t}+((\vel\cdot\na)\mathcal{B},\mathcal{W})_{\Omega_s^t}-((\na \vel)^T \mathcal{B},\mathcal{W})_{\Omega_s^t}- (\mathcal{B} \na \vel,\mathcal{W})_{\Omega_s^t}=0~.
\end{array}\right.\ee

Although the weak formulation \eqref{weakform} is straightforward, the time-dependent domains $\Omega_s^t$ and $\Omega_f^t$ introduce significant complexity in numerical discretization. To address this, we employ an ALE framework to map all variables back to a fixed reference domain $\hat{\Omega}$, simplifying the temporal discretization process.   We introduce the ALE mapping over the space-time domain:
\be
\mathcal{A}: \hat{\Omega}\times[0,T]\rightarrow \Omega~,\quad (\hat{\bx},t) \mapsto \mathcal{A}_t(\hat{\bx})~,
\ee
where the spatial transformation $\mathcal{A}_t(\hat{\bx}): \hat{\Omega} \rightarrow \Omega$ at any given time $t$ maps a reference structural domain $\hat{\Omega}_s$ to the current structural domain, satisfying $\mathcal{A}_t(\hat{\Omega}_s)=\Omega_s^t$ and $\mathcal{A}_t(\hat{\Omega}_f)=\Omega_f^t$, where $\hat{\Omega}_f = \hat{\Omega} \setminus \overline{\hat{\Omega}_s}$. When the reference domain is chosen as the initial configuration, $\hat{\Omega}_s = \hat{\Omega}_s$ and $\mathcal{A}_0=\mathbf{I}$. We also introduce the mesh velocity $w(\mathcal{A}(\hat{\bx},t),t)=\frac{\partial \mathcal{A}(\hat{\bx},t)}{\partial t}$ to construct the mapping through solving a harmonic extension problem
\begin{equation}
  \left\{\begin{array}{ll}
    \Delta w = 0~, & \mbox{in~}\Omega_f^t~,           \\[3mm]
    w = \vel~,     & \mbox{in~}\Omega_s^t~,           \\[3mm]
    w =0~,         & \mbox{on~}\partial\hat{\Omega}~.
  \end{array}\right.
\end{equation}

\begin{remark}
  In conventional FSI treatments, the structural domain $\Omega_s^t$ is often described strictly by a Lagrangian formulation, whereas the fluid domain $\Omega_f^t$ utilizes an ALE formulation. However, in our framework, the entire domain $\Omega$ is unified under a global ALE description. It is important to note that because the ALE mesh velocity ${w}$ is treated explicitly while the convective terms are treated implicitly in our temporal discretization, the relative convective velocity $(\mathbf{u} - {w})$ does not strictly vanish in $\Omega_s^t$ at the discrete level. Consequently, the structural domain is genuinely computed using an ALE formulation.
\end{remark}

\blue{
  By employing the ALE mapping, the time derivative of a physical quantity $f$ defined on the moving domain, evaluated with respect to a fixed reference coordinate $\hat{\bx}$, is defined as the ALE time derivative $\left.\frac{\partial f}{\partial t}\right|_{\hat{\bx}} = \frac{\partial f}{\partial t} + (w \cdot \na)f$. Consequently, the standard material derivative can be rewritten to explicitly account for the moving mesh by introducing the relative convective velocity $(\vel - w)$:
  \begin{equation}
    \frac{D f}{D t} = \frac{\partial f}{\partial t} + (\vel \cdot \na)f = \left.\frac{\partial f}{\partial t}\right|_{\hat{\bx}} + ((\vel - w) \cdot \na)f~.
  \end{equation}
  Substituting this relation into \eqref{weakform}, we obtain the equivalent weak formulation expressed in the current ALE configuration: find $(\mathbf{u},P^{s},P^{f},\mathcal{B})\in \mathbf{V}^t\times Q^t_s\times Q^t_f\times\mathbb{B}^t$ such that for all $(\bv,q^s,q^f,\mathcal{W})\in \mathbf{V}^t\times Q^t_s\times Q^t_f\times\mathbb{B}^t$
  \be\left\{\begin{array}{ll}
    \label{ale_current_weak}
    (\left.\frac{\partial \vel}{\partial t}\right|_{\hat{\bx}}, \bv)_{\hat{\Omega}} +(((\vel - w)\cdot \na)\vel,\bv)_{\hat{\Omega}} + \frac{1}{Re}(2\mathbf{D}(\mathbf{u}), \mathbf{D}(\bv))_{\Omega_f^t} \\
    ~~~~~~~~-(P^{s}, \na\cdot \bv)_{\Omega_s^t} -(P^{f}, \na\cdot \bv)_{\Omega_f^t} +\mu_s((\mathcal{B}-\mathbf{I}), \na\bv)_{\Omega_s^t}=0,                                                       \\[3mm]
    (\na\cdot \vel,q^s)_{\Omega_s^t}=0~,                                                                                                                                                                  \\[3mm]
    (\na\cdot \vel,q^f)_{\Omega_f^t}=0~,                                                                                                                                                                  \\[3mm]
    (\left.\frac{\partial \mathcal{B}}{\partial t}\right|_{\hat{\bx}},\mathcal{W})_{\Omega_s^t}+(((\vel - w)\cdot\na)\mathcal{B},\mathcal{W})_{\Omega_s^t}-((\na \vel)^T \mathcal{B},\mathcal{W})_{\Omega_s^t}- (\mathcal{B} \na \vel,\mathcal{W})_{\Omega_s^t}=0~.
  \end{array}\right.\ee
  This formulation clearly illustrates the convective effect of the moving mesh. In practical numerical implementations (as will be detailed in Remark 3.2), this current-configuration ALE formulation is exactly what is assembled and solved. However, to establish a rigorous temporal discretization scheme without the mathematical ambiguity of integrating over time-varying domains, it is necessary to formally pull all variables back to the fixed reference configuration $\hat{\Omega}$.
}

We denote $\hat{\vel}(\hat{\bx},t)=\vel(\mathcal{A}(\hat{\bx},t),t)$, $\hat{P}^s(\hat{\bx},t)=P^s(\mathcal{A}(\hat{\bx},t),t)$, $\hat{P}^f(\hat{\bx},t)=P^f(\mathcal{A}(\hat{\bx},t),t)$, $\hat{\mathcal{B}}(\hat{\bx},t)=\mathcal{B}(\mathcal{A}(\hat{\bx},t),t)$, $\hat{w}_{t}(\hat{\bx},t)=w_{t}(\mathcal{A}(\hat{\bx},t),t)$. Moreover, we introduce some notations in the reference configuration:
\be
\hat{J}=|\det(\frac{\partial \mathcal{A}_t}{\partial \hat{\bx}})|~,~\hat{\bF}=\frac{\partial \mathcal{A}_t}{\partial \hat{\bx}}~,~\hat{\bF}^{-1}=(\frac{\partial \mathcal{A}_t}{\partial \hat{\bx}})^{-1}~,\hat{\na}=\frac{\partial}{\partial \hat{\bx}}~, ~\hat{D}(\cdot) =\frac{1}{2}(\hat{\na}(\cdot)\hat{\bF}^{-1}+\hat{\bF}^{-T}(\hat{\na}(\cdot))^T)~.
\ee

Here, we define the corresponding function spaces on the reference domains as $\hat{\mathbf{V}}$, $\hat{Q}_s$, $\hat{Q}_f$, and $\hat{\mathbb{B}}$. Then the weak formulation \eqref{weakform} could be rewritten into the reference configuration as: find $(\hat{\mathbf{u}},\hat{P}^{s},\hat{P}^{f},\hat{\mathcal{B}},\hat{w})\in \hat{\mathbf{V}}\times \hat{Q}_s\times \hat{Q}_f\times\hat{\mathbb{B}}\times \hat{\mathbf{V}}$ such that for all $(\hat{\bv},\hat{q}^s,\hat{q}^f,\hat{\mathcal{W}},\hat{v})\in \hat{\mathbf{V}}\times \hat{Q}_s\times \hat{Q}_f\times\hat{\mathbb{B}}\times \hat{\mathbf{V}}$
\be\left\{\begin{array}{ll}
  \label{ref_weakform}
  (\hat{J}\frac{\partial \hat{\vel}}{\partial t}, \hat{\bv})_{\hat{\Omega}} +(\hat{J}(\hat{\vel} - \hat{w})\cdot \hat{\na}\hat{\vel}\hat{\bF}^{-1},\hat{\bv})_{\hat{\Omega}} + \frac{1}{Re}(\hat{J}2\hat{D}(\hat{\mathbf{u}})\hat{\bF}^{-T}, \hat{D}(\hat{\bv})\hat{\bF}^{-T})_{\hat{\Omega}_f}                                                                                                                      \\
  ~~~~~~~~-(\hat{J}\hat{P}^s, \operatorname{tr}(\hat{\na}\hat{\bv}\hat{\bF}^{-1}))_{\hat{\Omega}_s} -(\hat{J}\hat{P}^f, \operatorname{tr}(\hat{\na}\hat{\bv}\hat{\bF}^{-1}))_{\hat{\Omega}_f} +\mu_s(\hat{J}(\hat{\mathcal{B}}-\mathbf{I}), \hat{\na}(\hat{\bv})\hat{\bF}^{-1})_{\hat{\Omega}_s}=0~,                                                                                                          \\[3mm]
  (\hat{J}\operatorname{tr}(\hat{\na}\hat{\vel}\hat{\bF}^{-1}),\hat{q}^s)_{\hat{\Omega}_s}=0~,                                                                                                                                                                                                                                                                                                                       \\[3mm]
  (\hat{J}\operatorname{tr}(\hat{\na}\hat{\vel}\hat{\bF}^{-1}),\hat{q}^f)_{\hat{\Omega}_f}=0~,                                                                                                                                                                                                                                                                                                                       \\[3mm]
  (\hat{J}\frac{\partial \hat{\mathcal{B}}}{\partial t},\hat{\mathcal{W}})_{\hat{\Omega}_s}+(\hat{J}(\hat{\vel}-\hat{w})\cdot\hat{\na}\hat{\mathcal{B}}\hat{\bF}^{-1},\hat{\mathcal{W}})_{\hat{\Omega}_s}-(\hat{J}(\hat{\na} \hat{\vel}\hat{\bF}^{-1})^T \hat{\mathcal{B}},\hat{\mathcal{W}})_{\hat{\Omega}_s}-(\hat{J} \hat{\mathcal{B}} \hat{\na} \hat{\vel}\hat{\bF}^{-1},\hat{\mathcal{W}})_{\hat{\Omega}_s}=0~, \\
  (\hat{\na} \hat{w}, \hat{\na} \hat{v})_{\hat{\Omega}} = 0~,                                                                                                                                                                                                                                                                                                                                                        \\
  \hat{w}=\hat{\vel}~,~                                                                                           & \mbox{in~}\hat{\Omega}_s~,                                                                                                                                                                                                                                                                       \\
  \hat{w}=0~,~                                                                                                    & \mbox{on~}\partial\Omega~,                                                                                                                                                                                                                                                                       \\
  \frac{\partial \mathcal{A}(\hat{\bx},t)}{\partial t}=\hat{w}(\hat{\bx},t)~,~\mathcal{A}(\hat{\bx},0)=\hat{\bx}~ & \mbox{in~}\hat{\Omega}~.
\end{array}\right.\ee

\begin{remark}
  Pulling back the equations to a fixed reference domain $\hat{\Omega}$ conceptually transforms the moving-boundary problem into a fixed-domain system, thereby removing the ambiguity of differentiating in time over the moving domains $\Omega_s^t,\Omega_f^t$. This allows the coupled PDEs to be reduced to a standard system of ODEs after spatial discretization, providing a rigorous mathematical foundation for the design and analysis of high-order temporal schemes like implicit-explicit partitioned Runge-Kutta (IMEX-PRK) \cite{https://doi.org/10.48550/arxiv.2602.01094}. The reference configuration is thus used to \emph{define} the time-stepping scheme, not to assemble it in practice.
\end{remark}

\begin{remark}
  Once the temporal scheme has been defined on $\hat{\Omega}$, its implementation does not require solving in the reference domain. Since the ALE kinematic variables are treated explicitly, the current mesh at each time step is already known, and the change-of-variables formula relates the reference form \eqref{ref_weakform} to its current-configuration counterpart \eqref{ale_current_weak} at the continuous level. We therefore evaluate the finite element matrices directly on current coordinates, which avoids carrying the deformation Jacobians $\hat{\bF}^{-1}$ and $\hat{J}$ through every quadrature point. It should be emphasized that the mesh nodes move with the ALE map $\mathcal{A}_{h,t}$; no per-step remeshing or solution transfer between meshes is performed here, and remeshing is invoked only occasionally when mesh quality degrades.
\end{remark}

\section{Discretization}
Let $\mathcal{T}_h$ be a fitted triangulation of the reference domain $\Omega$ with mesh size $h$, which conforms to the reference fluid-structure interface $\hat{\Gamma} = \partial\hat{\Omega}_s \cap \partial\hat{\Omega}_f$. (Here, the reference configuration may be chosen as the initial configuration at $t=0$, or appropriately updated upon remeshing).

For spatial discretization, we adopt the classical Taylor-Hood $\mathcal{P}_2/\mathcal{P}_1$ finite element pair for velocity-pressure approximation, and $\mathcal{P}_1$ elements for the left Cauchy-Green tensor $\mathcal{B}$. Let $\hat{K}$ denote a reference simplex (triangle in 2D or tetrahedron in 3D), and let $\mathcal{P}_k(\hat{K})$ be the space of polynomials of total degree at most $k$ on $\hat{K}$.

For the discrete approximation of the ALE mapping, we use the finite element basis functions. Specifically, at each time $t$, $\mathcal{A}_{h,t} = \mathcal{A}_h(\cdot,t)$ is a piecewise polynomial map that moves the mesh nodes in a way that approximates the continuous motion. The discrete domain at time $t$ is then defined as
\[
  \Omega_h(t) = \mathcal{A}_{h,t}(\hat{\Omega}) = \Omega_{s,h}^t \cup \Omega_{f,h}^t \cup \Gamma_h^t~,
\]
where $\Omega_{s,h}^t = \mathcal{A}_{h,t}(\hat{\Omega}_s)$ and $\Omega_{f,h}^t = \mathcal{A}_{h,t}(\hat{\Omega}_f)$ are approximations to the true solid and fluid domains, respectively. The discrete interface $\Gamma_h^t = \partial\Omega_{s,h}^t \cap \partial\Omega_{f,h}^t$ approximates the true interface $\Gamma^t$. We also use the similar notations with subscript $h$ to denote the discrete ALE mapping related quantities, e.g., $\hat{\bF}_h = \frac{\partial \mathcal{A}_h}{\partial \hat{\bx}}$, $\hat{J}_h = |\det(\hat{\bF}_h)|$, etc. We shall use the same element for velocity to approximate the ALE mapping $\mathcal{A}_h$.

The finite element spaces on the reference domain $\hat{\Omega}$ are defined as:
\begin{align*}
  \hat{\mathbf{V}}_h & = \left\{ \hat{\mathbf{v}}_h \in {H}^1(\hat{\Omega})^d : \hat{\mathbf{v}}_h|_{K} \circ F_K^{-1} \in [\mathcal{P}_2(\hat{K})]^d, \ \forall K \in \mathcal{T}_h \right\},                                                 \\
  \hat{Q}_{s,h}      & = \left\{ \hat{q}_h^s \in L^2(\hat{\Omega}_s) : \hat{q}_h^s|_{K} \circ F_K^{-1} \in \mathcal{P}_1(\hat{K}), \ \forall K \in \mathcal{T}_h, K \subset \hat{\Omega}_s \right\},                                           \\
  \hat{Q}_{f,h}      & = \left\{ \hat{q}_h^f \in L^2(\hat{\Omega}_f) : \hat{q}_h^f|_{K} \circ F_K^{-1} \in \mathcal{P}_1(\hat{K}), \ \forall K \in \mathcal{T}_h, K \subset \hat{\Omega}_f \right\},                                           \\
  \hat{\mathbb{B}}_h & = \left\{ \hat{\mathbf{B}}_h \in [H^1(\hat{\Omega}_s)]^{d\times d} : \hat{\mathbf{B}}_h|_{K} \circ F_K^{-1} \in [\mathcal{P}_1(\hat{K})]^{d\times d}, \ \forall K \in \mathcal{T}_h, K \subset \hat{\Omega}_s \right\}, \\
\end{align*}
where $F_K: \hat{K} \to K$ is the isoparametric mapping from the reference element to each physical element $K$ in the reference mesh.
\begin{remark}
  By employing isoparametric elements, the second-order geometric approximation remains consistent with the finite element spaces, thereby minimizing boundary errors and preserving high-order spatial convergence for the FSI model.
\end{remark}

The space discretized problem is to find $(\hat{\mathbf{u}}_h,\hat{P}^{s}_h,\hat{P}^{f}_h,\hat{\mathcal{B}}_h,\hat{w}_{t,h})\in \hat{\mathbf{V}}_h\times \hat{Q}_{s,h}\times \hat{Q}_{f,h}\times \hat{\mathbb{B}}_h\times \hat{\mathbf{V}}_h$ such that for all $(\hat{\bv}_h,\hat{q}^s_h,\hat{q}^f_h,\hat{\mathcal{W}}_h,\hat{v}_h)\in \hat{\mathbf{V}}_h\times \hat{Q}_{s,h}\times \hat{Q}_{f,h}\times \hat{\mathbb{B}}_h\times \hat{\mathbf{V}}_h$
\be
\scalebox{0.8}{
  $\displaystyle
    \left\{\begin{array}{ll}
      \label{space_discrete_ref}
      (\hat{J}_h\frac{\partial \hat{\vel}_h}{\partial t}, \hat{\bv}_h)_{\hat{\Omega}} +(\hat{J}_h(\hat{\vel}_h - \hat{w}_{h})\cdot \hat{\na}\hat{\vel}_h\hat{\bF}_h^{-1},\hat{\bv}_h)_{\hat{\Omega}} + \frac{1}{Re}(\hat{J}_h2\hat{D}(\hat{\mathbf{u}}_h), \hat{D}(\hat{\bv}_h))_{\hat{\Omega}_f}                                                                                                                                                                \\
      ~~~~~~~~-(\hat{J}_h\hat{P}^s_h, \operatorname{tr}(\hat{\na}\hat{\bv}_h\hat{\bF}_h^{-1}))_{\hat{\Omega}_s} -(\hat{J}_h\hat{P}^f_h, \operatorname{tr}(\hat{\na}\hat{\bv}_h\hat{\bF}_h^{-1}))_{\hat{\Omega}_f} +\mu_s(\hat{J}_h(\hat{\mathcal{B}}_h-\mathbf{I}), \hat{\na}(\hat{\bv}_h)\hat{\bF}_h^{-1})_{\hat{\Omega}_s}=0~,                                                                                                                          \\[3mm]
      (\hat{J}_h\operatorname{tr}(\hat{\na}\hat{\vel}_h\hat{\bF}_h^{-1}),\hat{q}^s_h)_{\hat{\Omega}_s}=0~,                                                                                                                                                                                                                                                                                                                                                       \\[3mm]
      (\hat{J}_h\operatorname{tr}(\hat{\na}\hat{\vel}_h\hat{\bF}_h^{-1}),\hat{q}^f_h)_{\hat{\Omega}_f}=0~,                                                                                                                                                                                                                                                                                                                                                       \\[3mm]
      (\hat{J}_h\frac{\partial \hat{\mathcal{B}}_h}{\partial t},\hat{\mathcal{W}}_h)_{\hat{\Omega}_s}+(\hat{J}_h(\hat{\vel}_h-\hat{w}_{h})\cdot\hat{\na}\hat{\mathcal{B}}_h\hat{\bF}_h^{-1},\hat{\mathcal{W}}_h)_{\hat{\Omega}_s}-(\hat{J}_h(\hat{\na} \hat{\vel}_h\hat{\bF}_h^{-1})^T \hat{\mathcal{B}}_h,\hat{\mathcal{W}}_h)_{\hat{\Omega}_s}-(\hat{J}_h \hat{\mathcal{B}}_h \hat{\na} \hat{\vel}_h\hat{\bF}_h^{-1},\hat{\mathcal{W}}_h)_{\hat{\Omega}_s}=0~, \\
      (\hat{\na} \hat{w}_{h}, \hat{\na} \hat{v}_h)_{\Omega} = 0~,                                                                                                                                                                                                                                                                                                                                                                                                \\
      \hat{w}_{h}=\hat{\vel}_h~,~                                                                     & \mbox{in~}\hat{\Omega}_s~,                                                                                                                                                                                                                                                                                                                               \\
      \hat{w}_h = 0 ~,~                                                                               & \mbox{on~}\partial\Omega~,                                                                                                                                                                                                                                                                                                                               \\
      \frac{\partial \mathcal{A}_h(\bx,t)}{\partial t}=\hat{w}_{h}(\bx,t)~,~\mathcal{A}_h(\bx,0)=\bx~ & \mbox{in~}\Omega~,
    \end{array}\right.
  $
}
\ee
Now all of the variables are defined on the fixed reference domain, and we can perform the time discretization directly. We use a partitioned Runge-Kutta method for time discretization. For notational simplicity and algorithmic clarity, we recast the coupled system \eqref{space_discrete_ref} into the compact abstract form \eqref{star} by introducing the composite solution vector $\bX = (\hat{\vel}_h, \hat{P^s_h}, \hat{P^f_h}, \hat{\mathcal{B}}_h)$, test function vector $\Phi = (\hat{\bv}_h, \hat{q}^s_h, \hat{q}^f_h, \hat{\mathcal{W}}_h)$, the composite function space $\mathcal{X}_h = \hat{\mathbf{V}}_h \times \hat{Q}_{s,h} \times \hat{Q}_{f,h} \times \hat{\mathbb{B}}_h$, and the nonlinear operator $\mathcal{F}$ encapsulating the weak form of the governing equations:
\be
\label{star}
\begin{array}{ll}
  \mathcal{F}(\partial_t\bX_h,\bX_h; \Phi;\mathcal{A}_h, w_{h}) = 0~, & \mbox{for~all~}\Phi\in \mathcal{X}_h~, \\
  (\nabla w_{h}, \nabla v_h)_{\Omega} = 0~,                           & \mbox{for~all~}v_h\in \hat{V}_h~,      \\
  \hat{w}_{h} = \hat{\vel}_h~,                                        & \mbox{in~}\hat{\Omega}_s~,             \\
  \partial_t\mathcal{A}_h = \hat{w}_{h}~,                             & \mbox{in~}\hat{\Omega}~,               \\
  \mathcal{A}_h(\bx,0) = \bx~,                                        & \mbox{in~}\hat{\Omega}~,
\end{array}
\ee
where $\mathcal{F}$ is defined by
\be
\label{refer_discrete_F}
\begin{aligned}
  \mathcal{F}(\partial_t\bX_h,\bX_h; \Phi;\mathcal{A}_h, w_{h}) & =
  (\hat{J}_h\partial_t \hat{\vel}_h, \hat{\bv}_h)_{\hat{\Omega}} +(\hat{J}_h(\hat{\vel}_h - \hat{w}_{h})\cdot \hat{\na}\hat{\vel}_h\hat{\bF}_h^{-1},\hat{\bv}_h)_{\hat{\Omega}} + \frac{1}{Re}(\hat{J}_h2\hat{D}(\hat{\mathbf{u}}_h), \hat{D}(\hat{\bv}_h))_{\hat{\Omega}_f}                                                                                                             \\
                                                                & -(\hat{J}_h\hat{P}^s_h, \operatorname{tr}(\hat{\na}\hat{\bv}_h\hat{\bF}_h^{-1}))_{\hat{\Omega}_s} -(\hat{J}_h\hat{P}^f_h, \operatorname{tr}(\hat{\na}\hat{\bv}_h\hat{\bF}_h^{-1}))_{\hat{\Omega}_f} +\mu_s(\hat{J}_h(\hat{\mathcal{B}}_h-\mathbf{I}), \hat{\na}(\hat{\bv}_h)\hat{\bF}_h^{-1})_{\hat{\Omega}_s}~ \\
                                                                & +(\hat{J}_h\operatorname{tr}(\hat{\na}\hat{\vel}_h\hat{\bF}_h^{-1}),\hat{q}^s_h)_{\hat{\Omega}_s}+(\hat{J}_h\operatorname{tr}(\hat{\na}\hat{\vel}_h\hat{\bF}_h^{-1}),\hat{q}^f_h)_{\hat{\Omega}_f}~                                                                                                                    \\
                                                                & +(\hat{J}_h\partial_t \hat{\mathcal{B}}_h,\hat{\mathcal{W}}_h)_{\hat{\Omega}_s}+(\hat{J}_h(\hat{\vel}_h-\hat{w}_{h})\cdot\hat{\na}\hat{\mathcal{B}}_h\hat{\bF}_h^{-1},\hat{\mathcal{W}}_h)_{\hat{\Omega}_s}                                                                                                            \\
                                                                & -(\hat{J}_h(\hat{\na} \hat{\vel}_h\hat{\bF}_h^{-1})^T \hat{\mathcal{B}}_h,\hat{\mathcal{W}}_h)_{\hat{\Omega}_s}-(\hat{J}_h \hat{\mathcal{B}}_h \hat{\na} \hat{\vel}_h\hat{\bF}_h^{-1},\hat{\mathcal{W}}_h)_{\hat{\Omega}_s}~.
\end{aligned}
\ee

As discussed in the previous section, while the abstract operator $\mathcal{F}$ is formally constructed on the fixed reference domain $\hat{\Omega}$ to provide a rigorous basis for temporal discretization, it can be equivalently pushed forward to the current discrete spatial domain $\Omega_h^t = \mathcal{A}_h(\hat{\Omega}, t)$. By applying standard change-of-variables rules for integration, the geometric mapping terms $\hat{\bF}_h$ and $\hat{J}_h$ naturally vanish, and the nonlinear operator $\mathcal{F}$ evaluates identically to the weak form on the moving mesh:
\be
\label{current_discrete_F}
\begin{aligned}
  \mathcal{F}(\partial_t\bX_h,\bX_h; \Phi;\mathcal{A}_h, w_{h}) & =
  \left(\left.\frac{\partial \vel_h}{\partial t}\right|_{\hat{\bx}}, \bv_h\right)_{\Omega_h^t} +((\vel_h - w_{h})\cdot \na\vel_h,\bv_h)_{\Omega_h^t} + \frac{1}{Re}(2\mathbf{D}(\mathbf{u}_h), \mathbf{D}(\bv_h))_{\Omega_{f,h}^t}                     \\
                                                                & -(P^s_h, \na\cdot \bv_h)_{\Omega_{s,h}^t} -(P^f_h, \na\cdot \bv_h)_{\Omega_{f,h}^t} +\mu_s(\mathcal{B}_h-\mathbf{I}, \na\bv_h)_{\Omega_{s,h}^t}                               \\
                                                                & +(\na\cdot \vel_h,q^s_h)_{\Omega_{s,h}^t}+(\na\cdot \vel_h,q^f_h)_{\Omega_{f,h}^t}                                                                                                   \\
                                                                & +\left(\left.\frac{\partial \mathcal{B}_h}{\partial t}\right|_{\hat{\bx}},\mathcal{W}_h\right)_{\Omega_{s,h}^t}+((\vel_h-w_{h})\cdot\na\mathcal{B}_h,\mathcal{W}_h)_{\Omega_{s,h}^t} \\
                                                                & -((\na \vel_h)^T \mathcal{B}_h,\mathcal{W}_h)_{\Omega_{s,h}^t}-( \mathcal{B}_h \na \vel_h,\mathcal{W}_h)_{\Omega_{s,h}^t}~.
\end{aligned}
\ee
This equivalence is the cornerstone of the practical implementation. Evaluating the system residual does not require the explicit computation of the deformation gradient $\hat{\bF}_h$ or its determinant $\hat{J}_h$. The assembly proceeds with the standard spatial gradients ($\na$) and inner products over the updated mesh coordinates of $\Omega_h^t$ at each time step, bridging the reference-domain formulation with an efficient current-configuration implementation.

Similarly, we introduce an abstract operator $\mathcal{H}$ for the mesh velocity harmonic extension problem. Finding $w_h$ that satisfies the harmonic extension given a boundary/interface velocity $\mathbf{u}^*_h$ is denoted as $\mathcal{H}(w_h, \mathbf{u}^*_h; v_h) = 0$ for all $v_h \in \hat{\mathbf{V}}_h$, which encapsulates:

\vspace{-0.5em}
\be
\begin{array}{ll}
             (\nabla w_{h}, \nabla v_h)_{\hat{\Omega}} = 0~, & \mbox{for~all~}v_h\in \hat{\mathbf{V}}_h~, \\
             w_{h} = \mathbf{u}^*_h~,                        & \mbox{in~}\hat{\Omega}_s~,                              \\
             w_{h} = 0 ~,                                    & \mbox{on~}\partial\hat{\Omega}~.
\end{array}
\ee

\begin{remark}
  To ensure the exact algebraic equivalence of \eqref{ref_weakform} and \eqref{current_discrete_F}, the finite element spaces must remain invariant under the ALE mapping. A non-linear $\mathcal{A}_h$ would map the polynomial basis functions to non-polynomial functions, fundamentally distorting the trial and test spaces during the pull-back/push-forward operations. By enforcing $\mathcal{A}_h$ to be a piecewise linear (affine) mapping, we guarantee that the polynomial degree of the discrete spaces is strictly preserved across configurations. This affine constraint is maintained even for elements adjacent to curved interfaces to ensure the exact isomorphism of the functional spaces and the resulting discrete systems.
\end{remark}

We now proceed to discretize \eqref{star} in time using a semi-implicit Euler scheme. The idea is to set up a partitioned Runge-Kutta method \cite{ascher1997implicit,kennedy2003additive} where the mesh motion is solved separately from the main FSI system at each time step \cite{https://doi.org/10.48550/arxiv.2602.01094}.
Let $t_n = n\Delta t$ for $n=0,1,2,\ldots,N$ denote the discrete time levels with uniform step size $\Delta t$. Here and in what follows, the superscript $n$ in variables such as $\bX_h^n$ and $\mathcal{A}_h^n$ indicates the time step index. The time derivative is approximated as:
\be
\partial_t \bX_h^{n+1} \approx \frac{\bX_h^{n+1} - \bX_h^n}{\Delta t}~,~\partial_t \mathcal{A}_h^{n+1} \approx \frac{\mathcal{A}_h^{n+1} - \mathcal{A}_h^n}{\Delta t}~.
\ee

\begin{remark}
  We refer to the resulting scheme as a \emph{quasi-monolithic} formulation rather than a strictly monolithic one. The fluid momentum, the incompressibility constraints, the constitutive equation, and the evolution of $\mathcal{B}$ are assembled into a single nonlinear system $\mathcal{F}=0$ and solved implicitly by Newton's method, which preserves strong fluid-solid coupling and stability under matched-density conditions. The harmonic mesh extension $\mathcal{H}=0$ is treated explicitly through a staggered update at each Runge-Kutta stage. This decoupling of the mesh-motion sub-problem balances computational cost against time accuracy: when paired with the second-order IMEX-PRK scheme below, the global temporal order is preserved without forming the cross-coupled mesh-FSI Jacobian.
\end{remark}






Then the discretisation scheme for \eqref{star} reads:

\vspace{1em}
\begin{algorithm}[H]
  \caption{First-order semi-implicit Euler partitioned Runge-Kutta scheme}
  \For{$n=0, 1, \dots, N-1$}{
    Given $\bX_h^n$ and $\mathcal{A}_h^n$, find $\bX_h^{n+1}$ such that for all $\Phi\in \mathcal{X}_h$:
    \begin{equation}
      \label{time_discrete_new}
      \mathcal{F}\left(\frac{\bX_h^{n+1} - \bX_h^n}{\Delta t},\bX_h^{n+1}; \Phi;\mathcal{A}_h^n, w_{t,h}^n\right) = 0~.
    \end{equation}
    Find $w_{t,h}^{n+1}$ by solving the harmonic extension problem for all $v_h\in \hat{\mathbf{V}}_h$:
    \begin{equation}
      \label{time_discrete_wt_new}
      \mathcal{H}(w_{h}^{n+1}, \vel_h^{n+1}; v_h) = 0~.
    \end{equation}
    Update the mapping $\mathcal{A}_h^{n+1}$:
    \begin{equation}
      \label{time_discrete_A_new}
      \mathcal{A}_h^{n+1} = \mathcal{A}_h^n + \Delta t w_{h}^{n+1}~.
    \end{equation}
  }
\end{algorithm}
\vspace{1em}

We also show the second order time discretisation scheme as follows. Let $\gamma = 1-\frac{1}{\sqrt{2}}$, $\beta_0 = -\sqrt2$, $\beta_* = 1+\sqrt2$, $c_0 = -\frac{1}{\sqrt2}$, $c_* = 1+\frac{1}{\sqrt2}$, the second order time discretisation scheme reads:

\vspace{1em}
\begin{algorithm}[H]
  \caption{Second-order IMEX partitioned Runge-Kutta scheme}
  \For{$n=0, 1, \dots, N-1$}{
    Given $\bX_h^n$, find $w_{t,h}^{n*}$ by solving the harmonic extension problem for all $v_h\in \hat{\mathbf{V}}_h$:
    \begin{equation}
      \label{time_discrete_wt21_new}
      \mathcal{H}(w_{t,h}^{n*}, \vel_h^{n}; v_h) = 0~.
    \end{equation}
    Update the mapping $\mathcal{A}_h^{n*}$:\\
    \begin{equation}
      \label{time_discrete_A21_new}
      \mathcal{A}_h^{n*} = \mathcal{A}_h^n + \gamma\Delta t w_{h}^{n*}~.
    \end{equation}
    Find $\bX_h^{n*}$ such that for all $\Phi\in \mathcal{X}_h$:
    \begin{equation}
      \label{time_discrete21_new}
      \mathcal{F}\left(\frac{\bX_h^{n*} - \bX_h^n}{\gamma\Delta t},\bX_h^{n*}; \Phi;\mathcal{A}_h^{n*}, w_{h}^{n*}\right) = 0~.
    \end{equation}
    Find $w_{t,h}^{n+1}$ by solving the harmonic extension problem for all $v_h\in \hat{\mathbf{V}}_h$:
    \begin{equation}
      \label{time_discrete_wt22_new}
      \mathcal{H}(w_{h}^{n+1}, c_0\vel_h^{n} + c_* \vel_h^{n*}; v_h) = 0~.
    \end{equation}
    Update the mapping $\mathcal{A}_h^{n+1}$:
    \begin{equation}
      \label{time_discrete_A22_new}
      \mathcal{A}_h^{n+1} = \mathcal{A}_h^n + \Delta t w_{h}^{n+1}~.
    \end{equation}
    Find $\bX_h^{n+1}$ such that for all $\Phi\in \mathcal{X}_h$:
    \begin{equation}
      \label{time_discrete22_new}
      \mathcal{F}\left(\frac{\bX_h^{n+1} - \beta_0\bX_h^n-\beta_*\bX_h^{n*}}{\gamma\Delta t},\bX_h^{n+1}; \Phi;\mathcal{A}^{n+1}, w_{h}^{n+1}\right) = 0~.
    \end{equation}
  }
\end{algorithm}
\vspace{1em}

\section{Local Updating Algorithm for Multi-scale Simulation}
Microfluidic devices impose a scale separation in which the channel length spans tens to hundreds of times the particle diameter, so a global FSI computation on the full domain is prohibitive in three dimensions. The flow regime, however, is favorable. At low Reynolds number, the geometric confinement of the channel walls produces a hydrodynamic screening effect \cite{Blake1971} that decays the particle-induced velocity perturbation exponentially along the channel axis \cite{Liron1978}; beyond roughly 10 to 15 particle diameters the field returns to the undisturbed background laminar flow. We exploit this localization by confining the FSI solve, the moving mesh, and the deformation history $\mathcal{B}$ to a body-fitted sub-domain that follows the particle, and reading the surrounding flow from a precomputed steady background. The algorithm proceeds as follows.

\begin{algorithm}[H]
  \caption{Localized Updating Strategy for Multi-scale FSI}
  \label{alg:local_updating}
  Precompute the steady background velocity field $\mathbf{u}_{bg}$ on the static background mesh $\mathcal{M}_{bg}$.\\
  Generate an initial local body-fitted mesh $\mathcal{M}_{local}$ around the solid sphere.\\
  Initialize the coupled solution vector $\bX_{local}^0 = (\mathbf{u}, P^s, P^f, \mathcal{B})^0$ and mesh velocity $w_{t}^0$.\\
  \vspace{0.5em}
  \For{$n=0, 1, \dots, N-1$}{
    Interpolate the background velocity $\mathbf{u}_{bg}$ onto the boundary $\partial \mathcal{M}_{local}$ to enforce the Dirichlet boundary condition $\mathbf{u}_{\Gamma}$.\\
    Given $\bX_{local}^n$ on the current mesh $\mathcal{M}_{local}$, solve the quasi-monolithic FSI system with Newton's method to obtain the updated state $\bX_{local}^{n+1}$.\\
    Compute the mesh motion velocity $w_t^{n+1}$ based on $\mathbf{u}_{local}^{n+1}$.\\
    Update the sphere's position and the nodal coordinates of the mesh:
    \begin{equation}
      \label{update_mesh_coords}
      x_{local}^{n+1} = x_{local}^n + \Delta t w^{n+1}~.
    \end{equation}
    \eIf{the sphere displacement and mesh distortion exceed the predefined threshold}{
      Generate a new body-fitted local mesh $\mathcal{M}_{local}^*$ centered at the updated sphere position.\\
      Interpolate variables $(\mathbf{u}, \mathcal{B}, w_t)^{n+1}$ from the old mesh $\mathcal{M}_{local}$ to the newly generated mesh $\mathcal{M}_{local}^*$ within overlapping regions.\\
      In non-overlapping regions, populate the velocity field directly with $\mathbf{u}_{bg}$ and set other variables to zero.\\
      Replace the active mesh topology: $\mathcal{M}_{local} \leftarrow \mathcal{M}_{local}^*$.\\
    }{
      Retain the current mesh topology $\mathcal{M}_{local}$ (only nodal coordinates are updated).\\
    }
  }
\end{algorithm}

\begin{remark}
  The localized updating strategy reduces the global computational cost by confining resolution to the body-fitted sub-domain and reading the surrounding flow from the precomputed steady background.
  Field interpolation across mesh transitions preserves continuity of the velocity and the deformation history $\mathcal{B}$, sustaining numerical stability.
  The strategy is demonstrated for a 2D rigid sphere in a long channel and extends to 3D with appropriate mesh adaptation.
\end{remark}

\begin{figure}
  \centering
  \includegraphics[width=0.6\textwidth]{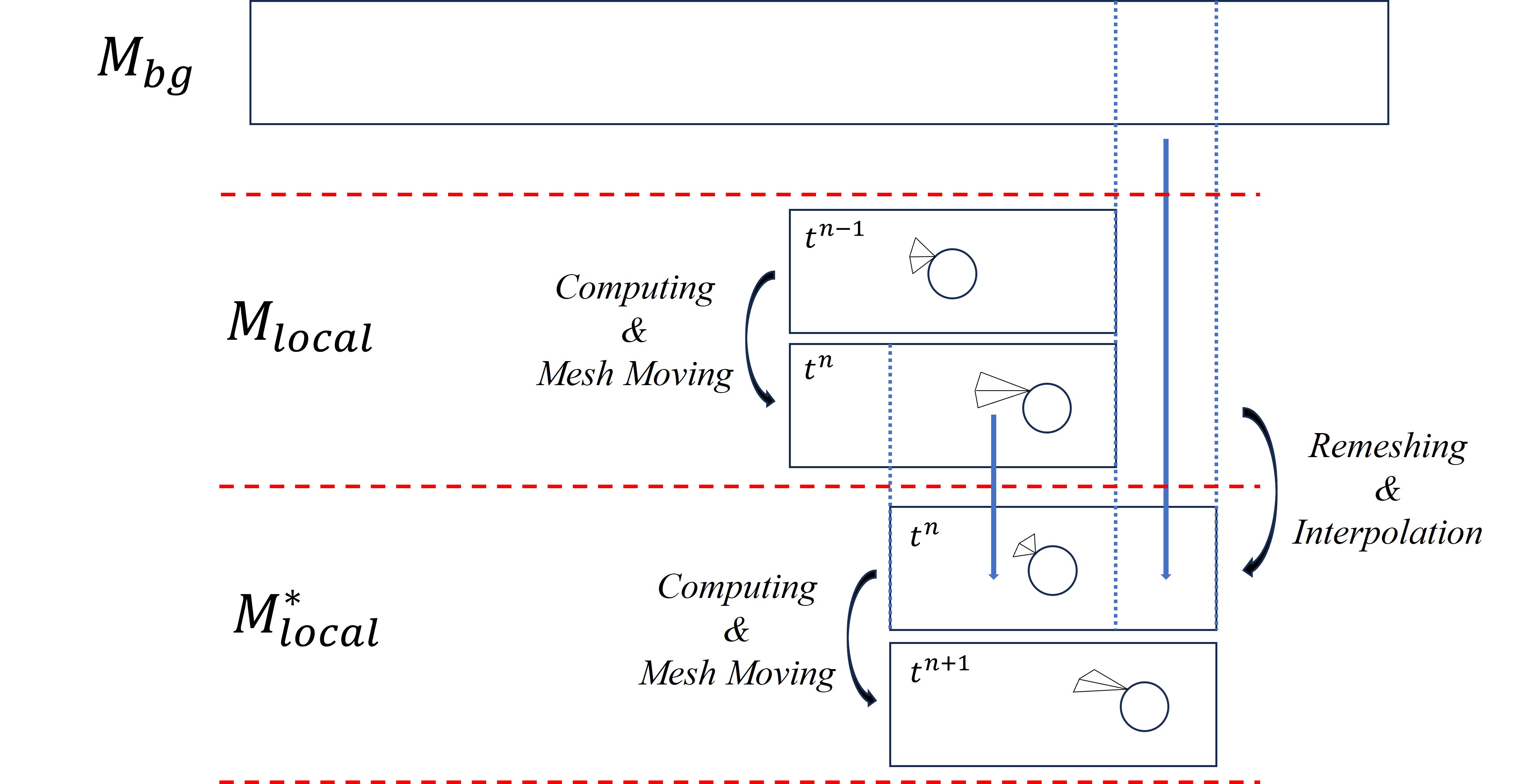}
  \caption{The figure illustrates the dynamic evolution of local mesh configurations across time steps, as well as the process of generating new meshes and interpolating field variables from the outdated mesh to the newly refined mesh when mesh quality metrics exceed predefined thresholds, thereby enabling continuous computational progression}
  \label{fig:algorithm}
\end{figure}

\begin{figure}
  \centering
  \includegraphics[width=0.6\textwidth]{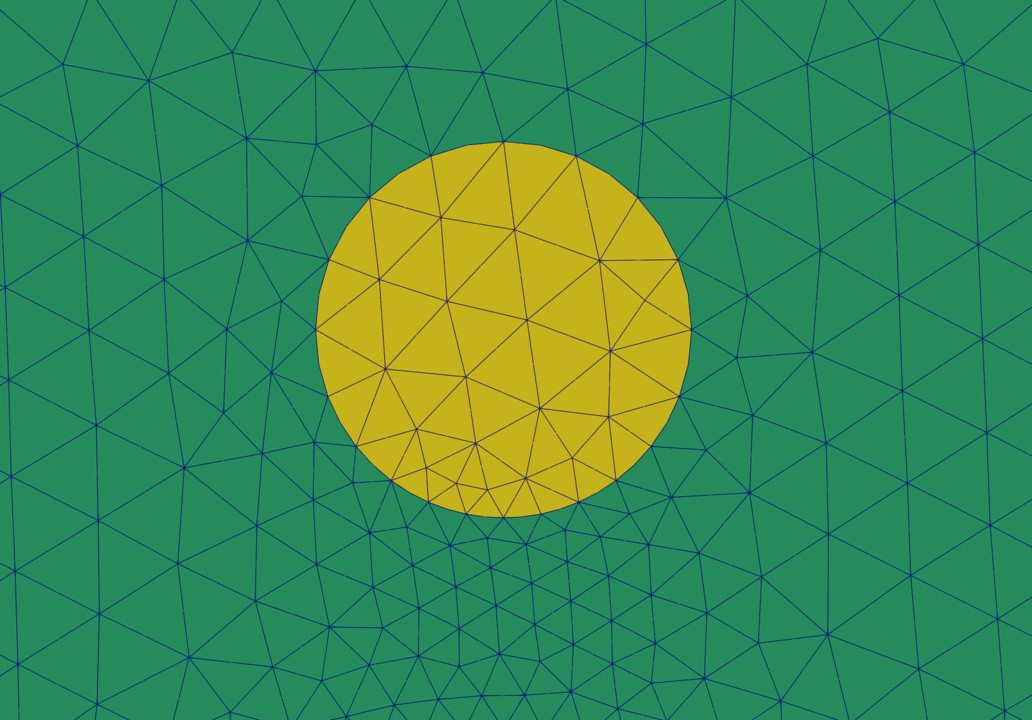}
  \caption{The second-ordered mesh fitting the boundary of the solid. The green elements represent the fluid phase while the yellow ones stand for solid phase. Notably, due to the application of high-order meshes, the elements at the fluid-solid interface adopt curved edges rather than straight edges. Unlike visualization representations, during finite element integration procedures, these curved edges are approximated using high-order polynomials, thereby better capturing curvature-induced geometric characteristics and achieving enhanced convergence rates in velocity fields.}
  \label{fig:2nd_ordered_mesh}
\end{figure}

Figure \ref{fig:algorithm} demonstrates how the velocity field is updated and how the mesh changes throughout the algorithm. Also, Second-order elements were employed to enhance the geometric fidelity in approximating curved boundaries of solid structures, as depicted in Figure \ref{fig:2nd_ordered_mesh}.

\section{Numerical Results}
This section validates the proposed framework on four problems of increasing complexity: a 2D convergence study with a rigid particle through two pillars, a 3D free-falling sphere benchmarked against the ten Cate experiment, the Turek-Hron FSI3 benchmark for matched-density added-mass coupling, and a multi-scale particle-focusing simulation in a spiral microchannel.

\subsection{Convergence study by double-pillar channel flow}
\begin{figure}
  \centering
  \includegraphics[width=0.6\textwidth]{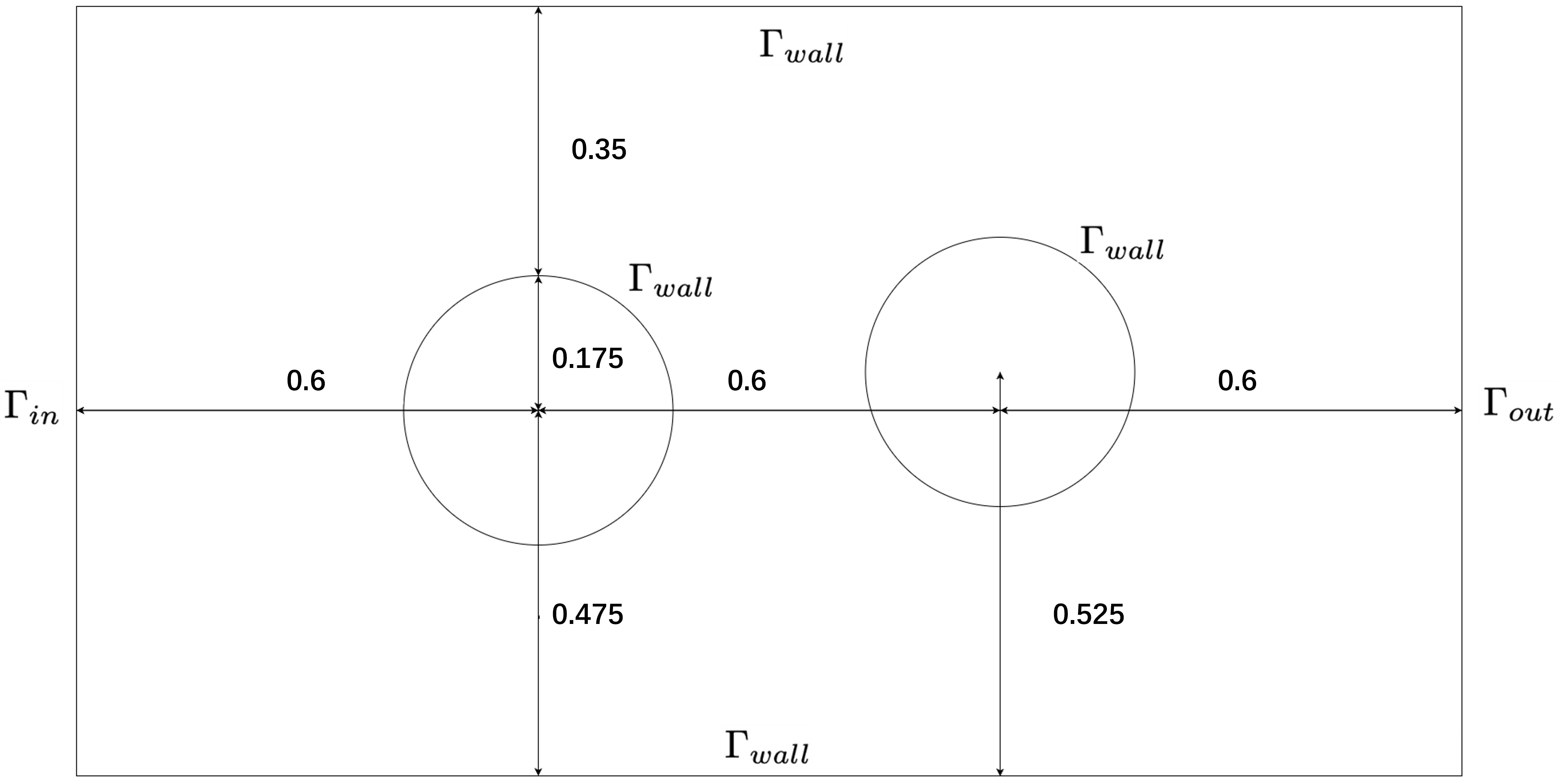}
  \caption{The domain set up for convergence study.}
  \label{fig:benchmark1}
\end{figure}

As illustrated in Figure \ref{fig:benchmark1}, we first design the following test case to verify the convergence of our algorithm. The calculation domain is set to be rectangular with width $W=1$ and length $L=1.8$. Inside the domain, two disk regions are cut to represent the pillars.

A spherical particle with radius $r=0.08$ is released from rest at (0.60, 0.76), while a parabolic inflow is imposed from left to right. The inflow boundary condition is set as

\be\label{eq:inflow-B1}
\mathbf{u} = 0.5u_0\frac{y(W-y)}{W^2}\mathbf{e}_x~.
\ee

The corresponding parameter values are summarized below:
\be
Re = 3~, ~\mu_s = 10^9~, ~u_0 = 8~.
\ee

\begin{figure}
  \centering
  \includegraphics[width=0.6\textwidth]{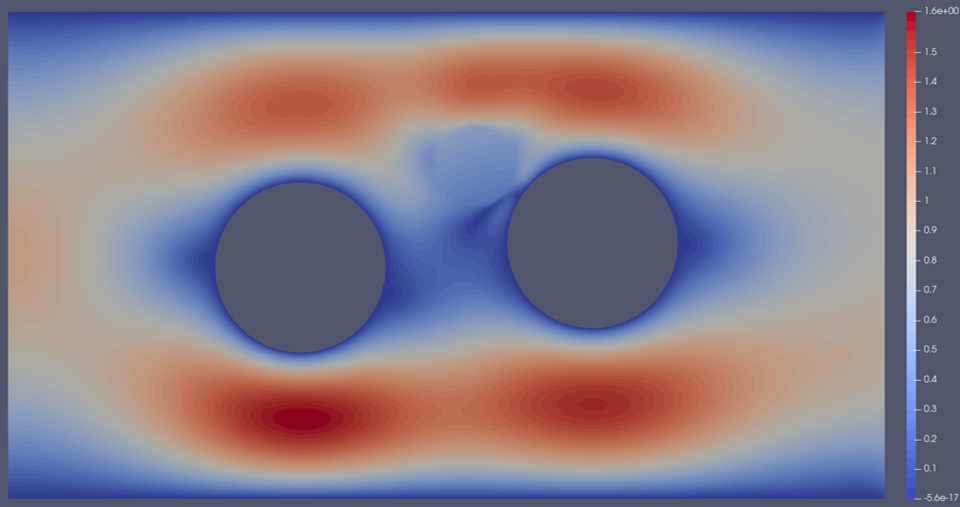}
  \caption{Magnitude of the velocity field obtained from the FSI simulation. The high-modulus particle perturbs the flow near the pillars when the gap is small.}
  \label{fig:benchmark1_velocity}
\end{figure}

Here, we employ a larger modulus to make the particle motion closer to that of a rigid body, thereby reducing the factors that may affect the order of convergence. Figure \ref{fig:benchmark1_velocity} presents the magnitude of the velocity field at a representative time instant in the simulation.

\begin{figure}
  \centering
  \includegraphics[width=1.0\textwidth]{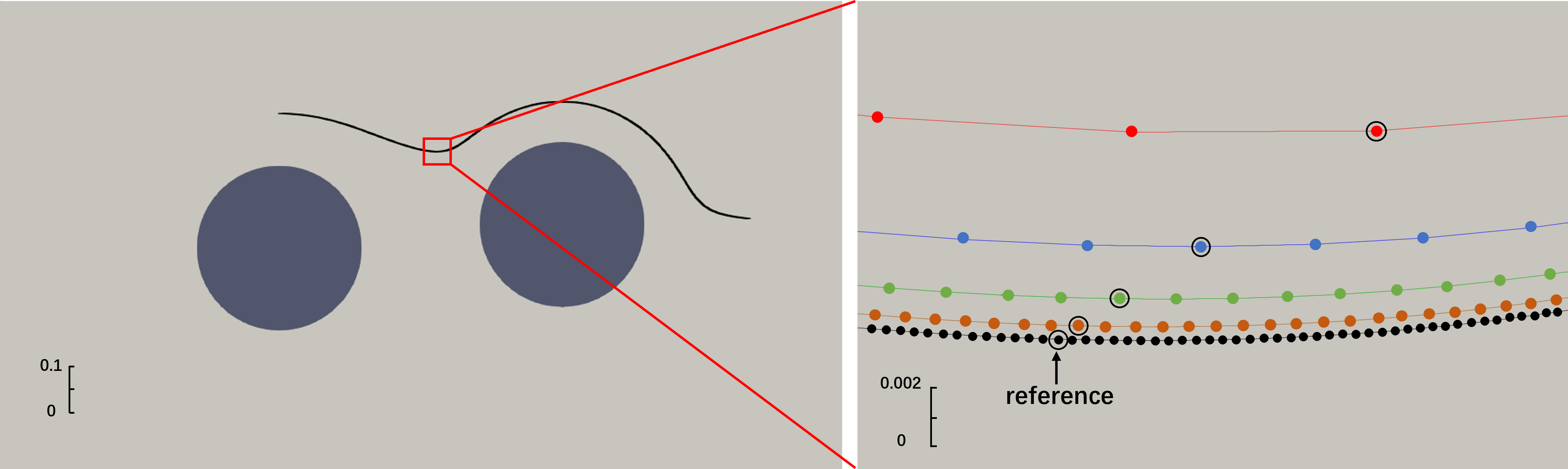}
  \caption{{The figure presents the numerically computed trajectory of the center of the rigid particle. The enlarged view shows the detail of the trajectories obtained with a mesh size of 4/100 and time-step sizes of 3/200 (red), 3/400 (blue), 3/800 (green), 3/1600 (brown) and 3/3200 (black) respectively. The particle positions at different time steps are marked by points along the curves. In particular, the point enclosed by the circle indicates the particle’s position at time t=0.375. For the convergence study, the particle position at this time is used to compute the numerical errors, where the solution obtained with the smallest time step ($\Delta t = 3/3200$) is taken as the reference.}}
  \label{fig:benchmark1_trajectory}
\end{figure}

For a prescribed total simulation time, we record the trajectory of the particle’s centroid throughout the whole process, as shown in Figure \ref{fig:benchmark1_trajectory}.

{The trajectories vary visibly with the mesh size and the time-step size, with the largest discrepancy occurring as the particle passes between the two pillars (Figure \ref{fig:benchmark1_trajectory}). Because the trajectory difference is dominated by the vertical component, we use the vertical coordinate of the centroid to evaluate the error:}

\begin{equation}
  \mathcal{E}(h, \Delta t) = \max_{t \in [0,T]} \, \big| \mathbf{y}_h(t) - \mathbf{y}_{\mathrm{ref}}(t) \big|~,
\end{equation}
where $\mathbf{y}_h(t)$ is the vertical centroid coordinate of the simulated particle and $\mathbf{y}_{\mathrm{ref}}(t)$ is the corresponding reference value, obtained by linear interpolation in time of the nodal values along the reference trajectory.


We first examined the convergence rates of the first-order backward Euler scheme and the second-order IMEX partitioned Runge--Kutta scheme described in Section 3. The results are shown in Tables \ref{table_3} and \ref{table_4}. The reference for each case is the trajectory computed with $\Delta t = 3/3200$.

\begin{table}[h]
  \centering
  \caption{Time convergence study of the first-order time integration scheme for a rigid particle in double-pillar channel flow ($Re = 3$, $\mu_s = 10^9$). The spatial mesh size is fixed at $h = 4/100$, and the reference trajectory is computed using a fine time step of $\Delta t = 3/3200$.}
  \label{table_3}
  \begin{tabular}{c c c}
    \toprule
    time step ($\Delta t$) & Error ($\mathcal{E}$) & convergence rate \\
    \midrule
    3/200                  & 6.9e-2                &                  \\
    3/400                  & 3.1e-2                & 1.15             \\
    3/800                  & 1.3e-2                & 1.25             \\
    3/1600                 & 4.4e-3                & 1.56             \\
    \bottomrule
  \end{tabular}
\end{table}

\begin{table}[h]
  \centering
  \caption{Time convergence study of the second-order time integration scheme for a rigid particle in double-pillar channel flow ($Re = 3$, $\mu_s = 10^9$). The spatial mesh size is fixed at $h = 4/100$, and the reference trajectory is computed using a fine time step of $\Delta t = 3/3200$.}
  \label{table_4}
  \begin{tabular}{c c c}
    \toprule
    time step ($\Delta t$) & Error ($\mathcal{E}$) & convergence rate \\
    \midrule
    3/200                  & 4.5e-2                &                  \\
    3/400                  & 1.1e-2                & 2.03             \\
    3/800                  & 2.8e-3                & 1.97             \\
    3/1600                 & 6.1e-4                & 2.20             \\
    \bottomrule
  \end{tabular}
\end{table}

The measured convergence rates match the analytical orders: the second-order IMEX-PRK scheme attains $\mathcal{O}(\Delta t^2)$, while the first-order backward Euler discretization is limited to $\mathcal{O}(\Delta t)$. The IMEX-PRK scheme delivers an order-of-magnitude reduction in trajectory error at every refinement level compared to backward Euler, which is consistent with the smaller numerical dissipation of the high-order scheme during the strong-coupling phase as the particle passes between the pillars.

Setting the time step to $\Delta t = 3/800$, we compare the spatial convergence of first- and second-order meshes to demonstrate the advantage of isoparametric finite elements. The corresponding results are summarized in Tables \ref{table_1} and \ref{table_2}. The reference trajectory for each case is computed with mesh size $h = 1/100$.

\begin{table}[h]
  \centering
  \caption{Spatial convergence study using first-order mesh. Time integration is performed using the first-order scheme with a fixed time step $\Delta t = 3/800$ ($Re = 3, \mu_s = 10^9$). The reference trajectory is calculated with $h = 1/100$.}
  \label{table_1}
  \begin{tabular}{c c c}
    \toprule
    Mesh size ($h$) & Error ($\mathcal{E}$) & Convergence rate \\
    \midrule
    4/100           & 1.6e-3                &                  \\
    $2\sqrt{2}/100$ & 8.3e-4                & 1.89             \\
    2/100           & 4.0e-4                & 2.11             \\
    $\sqrt{2}/100$  & 1.7e-4                & 2.47             \\
    \bottomrule
  \end{tabular}
\end{table}

\begin{table}[h]
  \centering
  \caption{Spatial convergence study using second-order mesh. Time integration is performed using the first-order scheme with a fixed time step $\Delta t = 3/800$ ($Re = 3, \mu_s = 10^9$). The reference trajectory is calculated with $h = 1/100$.}
  \label{table_2}
  \begin{tabular}{c c c}
    \toprule
    Mesh size ($h$) & Error ($\mathcal{E}$) & Convergence rate \\
    \midrule
    4/100           & 4.9e-4                &                  \\
    $2\sqrt{2}/100$ & 2.2e-4                & 2.31             \\
    2/100           & 8.1e-5                & 2.88             \\
    $\sqrt{2}/100$  & 2.3e-5                & 3.63             \\
    \bottomrule
  \end{tabular}
\end{table}

The second-order mesh raises the spatial convergence rate beyond the linear-element bound and lowers the trajectory error by nearly an order of magnitude at the finest resolution. This confirms that isoparametric high-order discretization preserves geometric fidelity at smooth interfaces under large structural translations, removing the staircase artifacts associated with linear or piecewise-constant interface representations.
The convergence study and the trajectory match establish the spatial and temporal accuracy of the method on the moving-interface FSI problem.




\subsection{3D FSI simulation of a free falling sphere in fluid}
Beyond demonstrating satisfactory performance in two-dimensional settings, the present model and its associated algorithm can be readily extended to three-dimensional cases.

  {Following the data provided in \cite{ten2002particle}, we designed and performed a simulation of the free fall of a small sphere in a tank. The tank has dimensions of $100\,\mathrm{mm}\times100\,\mathrm{mm}\times160\,\mathrm{mm}$. The sphere density is $1120\,\mathrm{kg/m^3}$, released from rest at point $(50\,\mathrm{mm}, 50\,\mathrm{mm}, 120\,\mathrm{mm})$. Simulations were performed for four working liquids of different viscosities; the corresponding density, dynamic viscosity, and Reynolds number for each case are summarized in Table~\ref{table_parameter}. A no-slip boundary condition is applied on all the walls of the tank. The dimensionless shear modulus of the ball is set to $\mu_s = 10^{9}$ to approximate rigid-body motion.}

\begin{table}[h]
  \centering
  \caption{Parameters of the four working liquids used in the free-fall benchmark, taken from \cite{ten2002particle}.} 
  \label{table_parameter}
  \begin{tabular}{c c c c} 
    \toprule
    Case & Density ($\mathrm{kg/m^3}$) & Dynamic viscosity ($\mathrm{mPa\cdot s}$) & Reynolds number \\
    \midrule
    1    & 970                  & 373                & 1.5             \\
    2    & 965                  & 212                & 4.1             \\
    3    & 962                  & 113                & 11.6            \\
    4    & 960                  & 58                 & 31.9            \\
    \bottomrule
  \end{tabular}
\end{table}

\begin{figure}[h]
  \centering
  \begin{subfigure}{0.23\textwidth}
    \includegraphics[width=\linewidth]{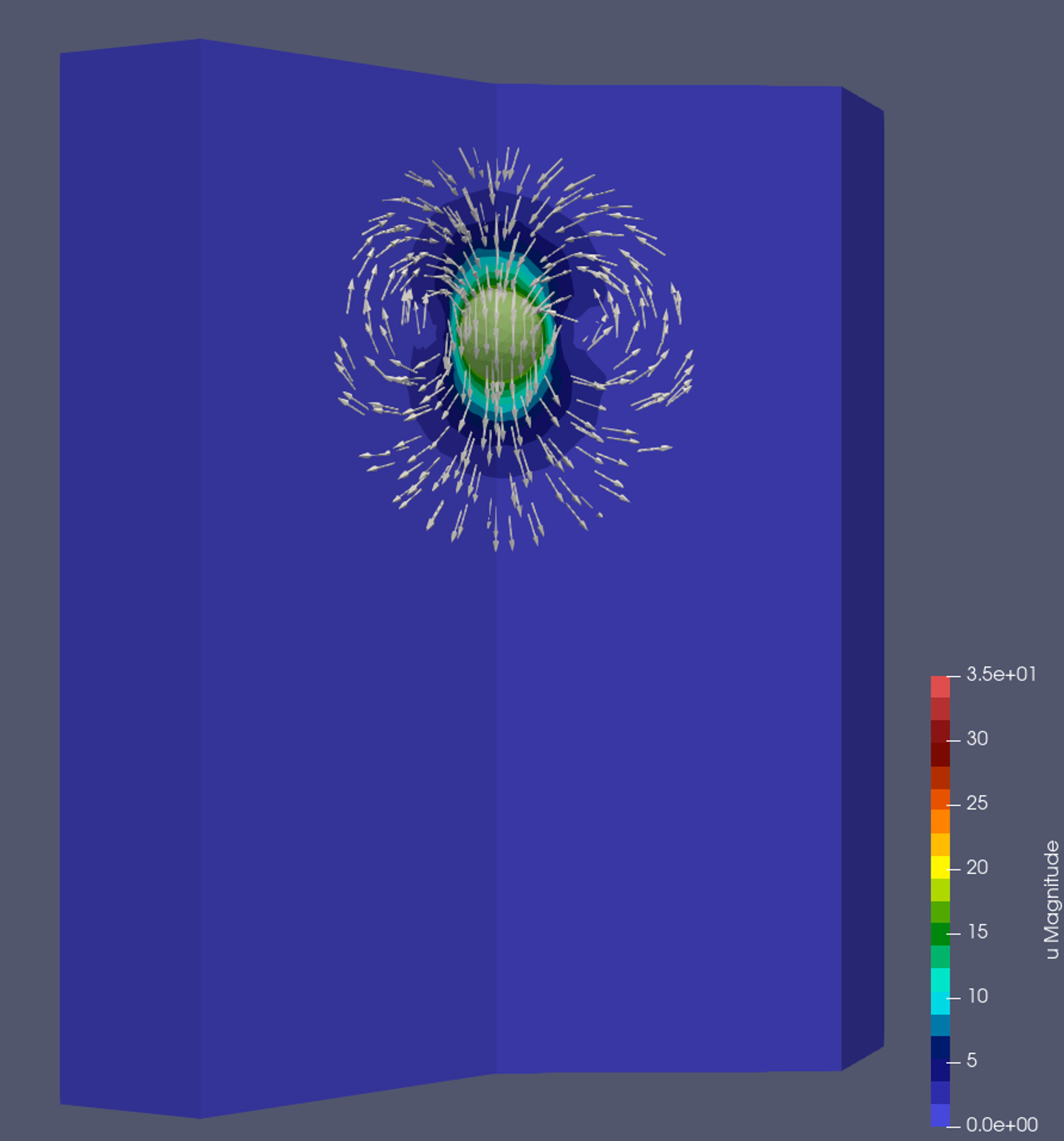}
    \caption{t=0.06s}
  \end{subfigure}\hfill
  \begin{subfigure}{0.23\textwidth}
    \includegraphics[width=\linewidth]{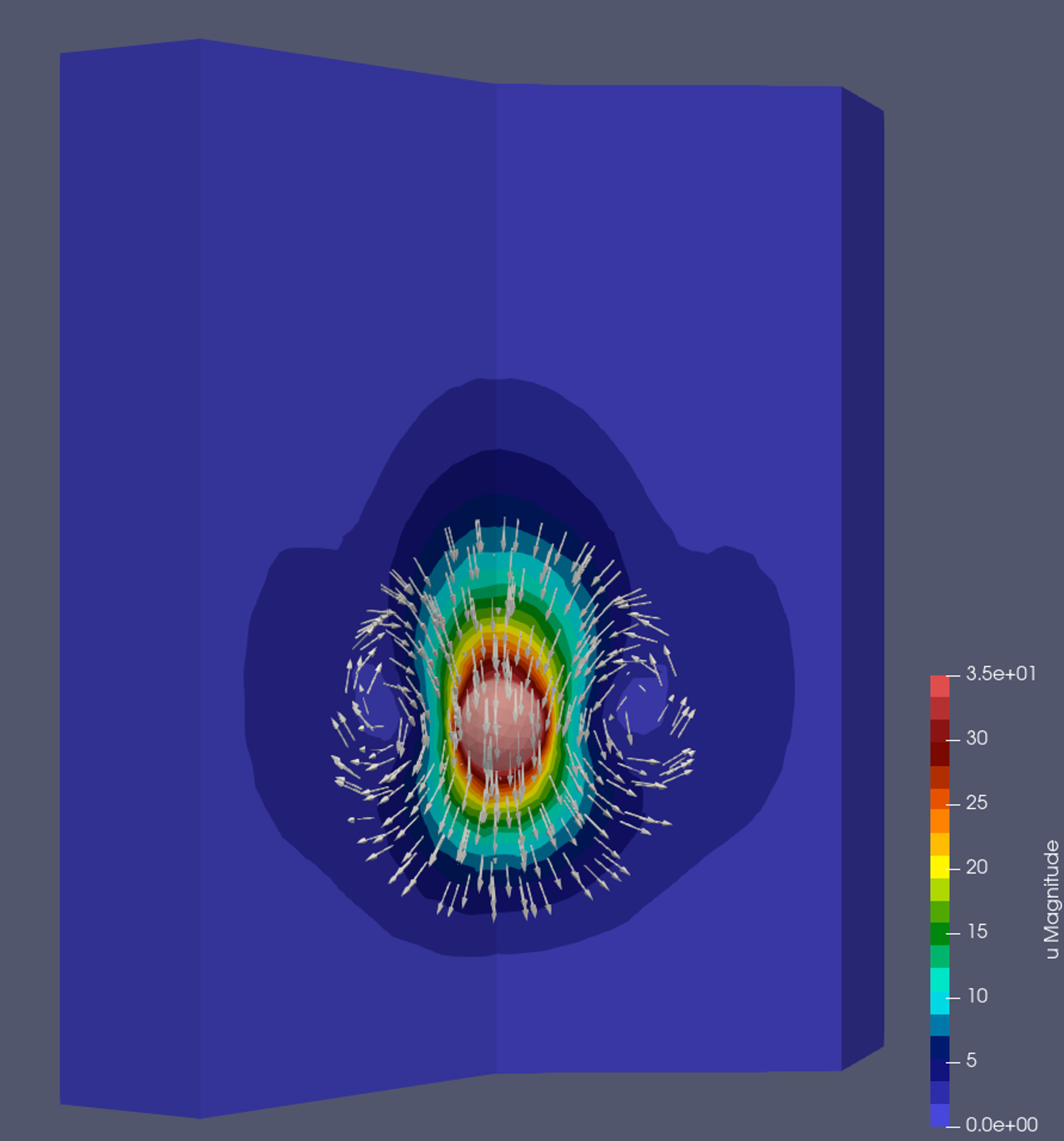}
    \caption{t=2s}
  \end{subfigure}
  \begin{subfigure}{0.23\textwidth}
    \includegraphics[width=\linewidth]{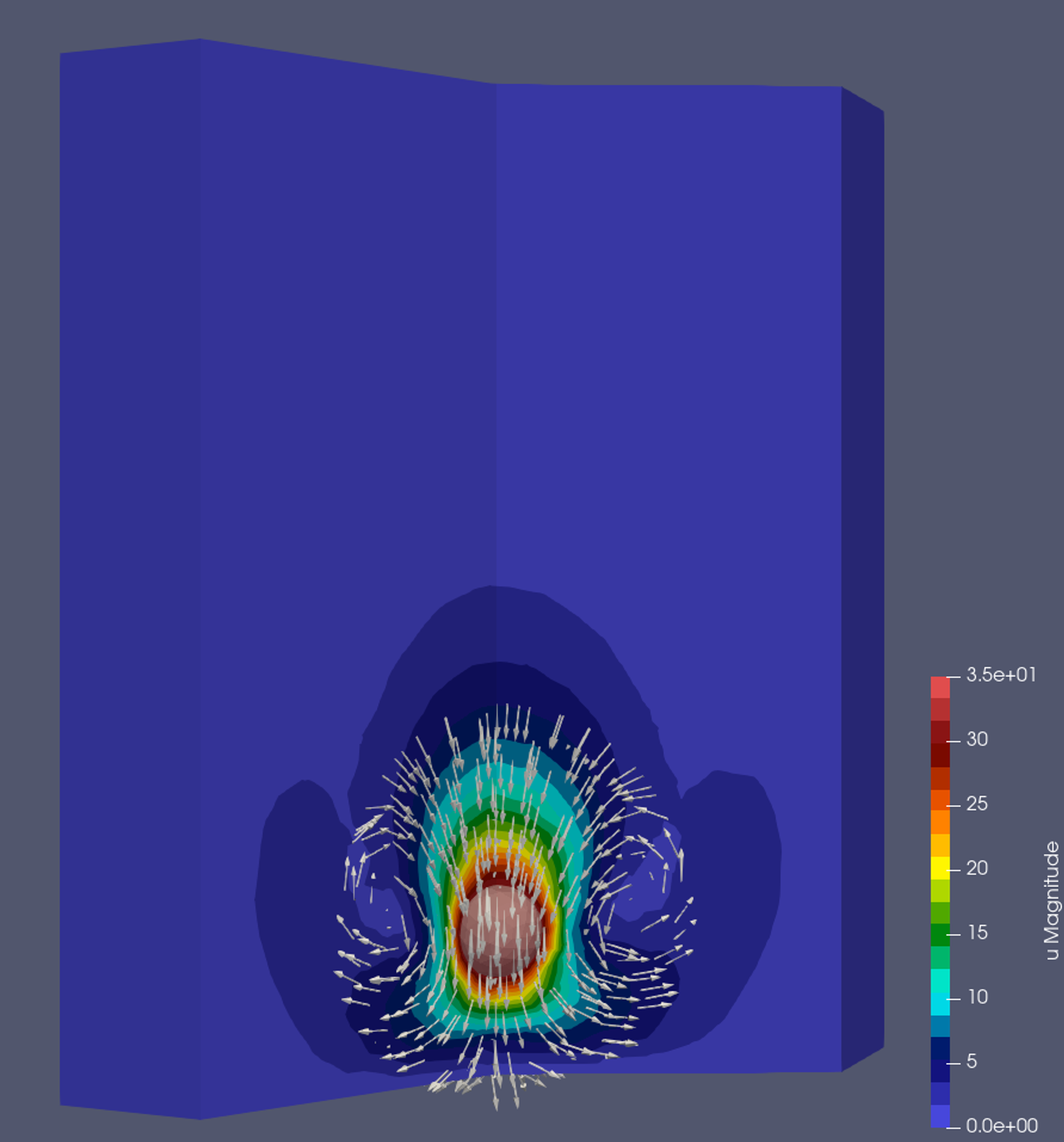}
    \caption{t=3s}
  \end{subfigure}\hfill
  \begin{subfigure}{0.23\textwidth}
    \includegraphics[width=\linewidth]{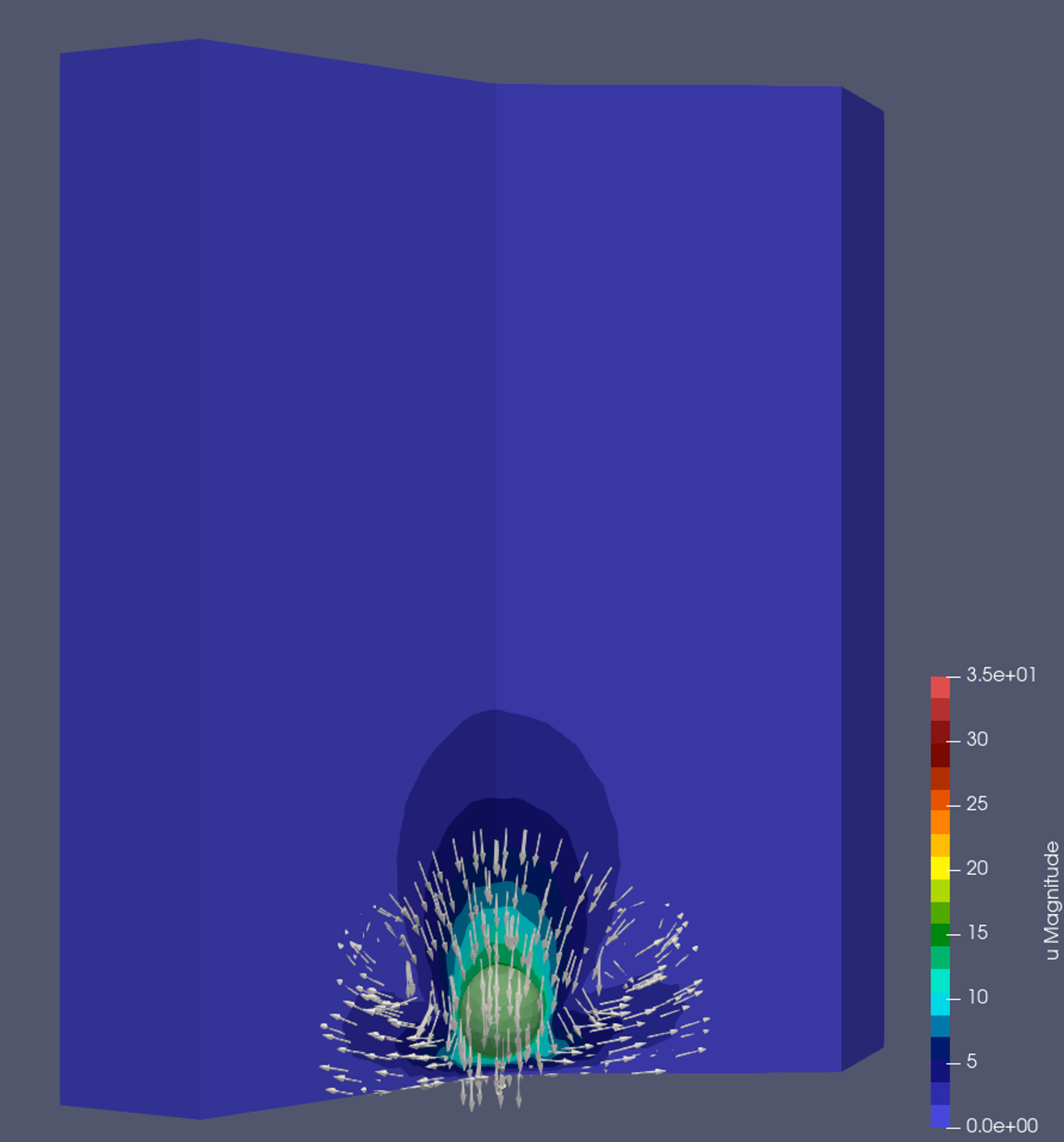}
    \caption{t=3.5s}
  \end{subfigure}

  \caption{The magnitude of the velocity field on the cross-section of the flow channel with different flux (left) along with the zoomed views (right) of the local part in the red frame.}
  \label{fig:freefall_different_time}
\end{figure}

\begin{figure}[htbp]
  \centering
  \begin{subfigure}{0.48\textwidth}
    \includegraphics[width=\linewidth]{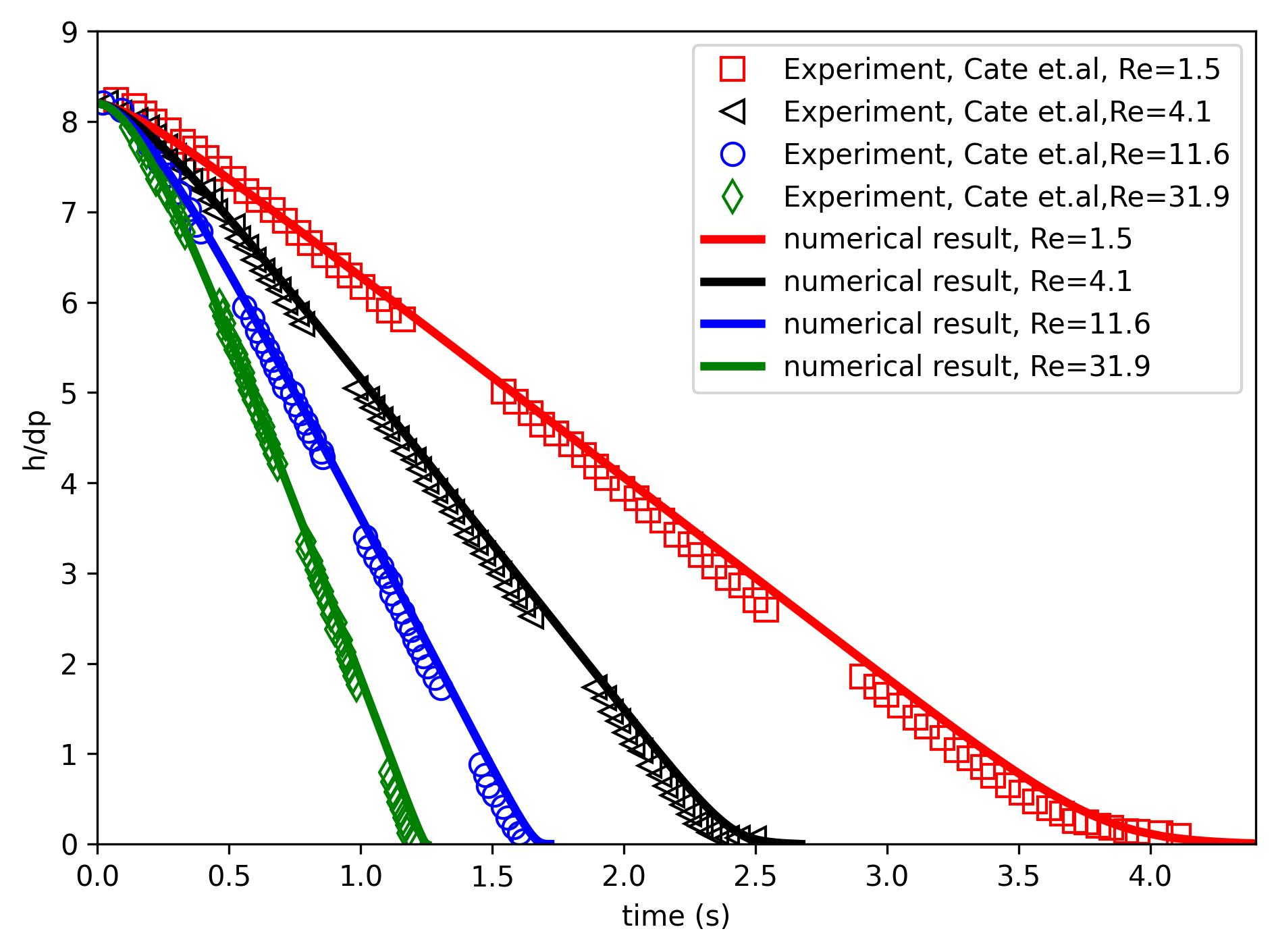}
    \caption{Height-diameter ratio vs time}
  \end{subfigure}
  \hfill
  \begin{subfigure}{0.48\textwidth}
    \includegraphics[width=\linewidth]{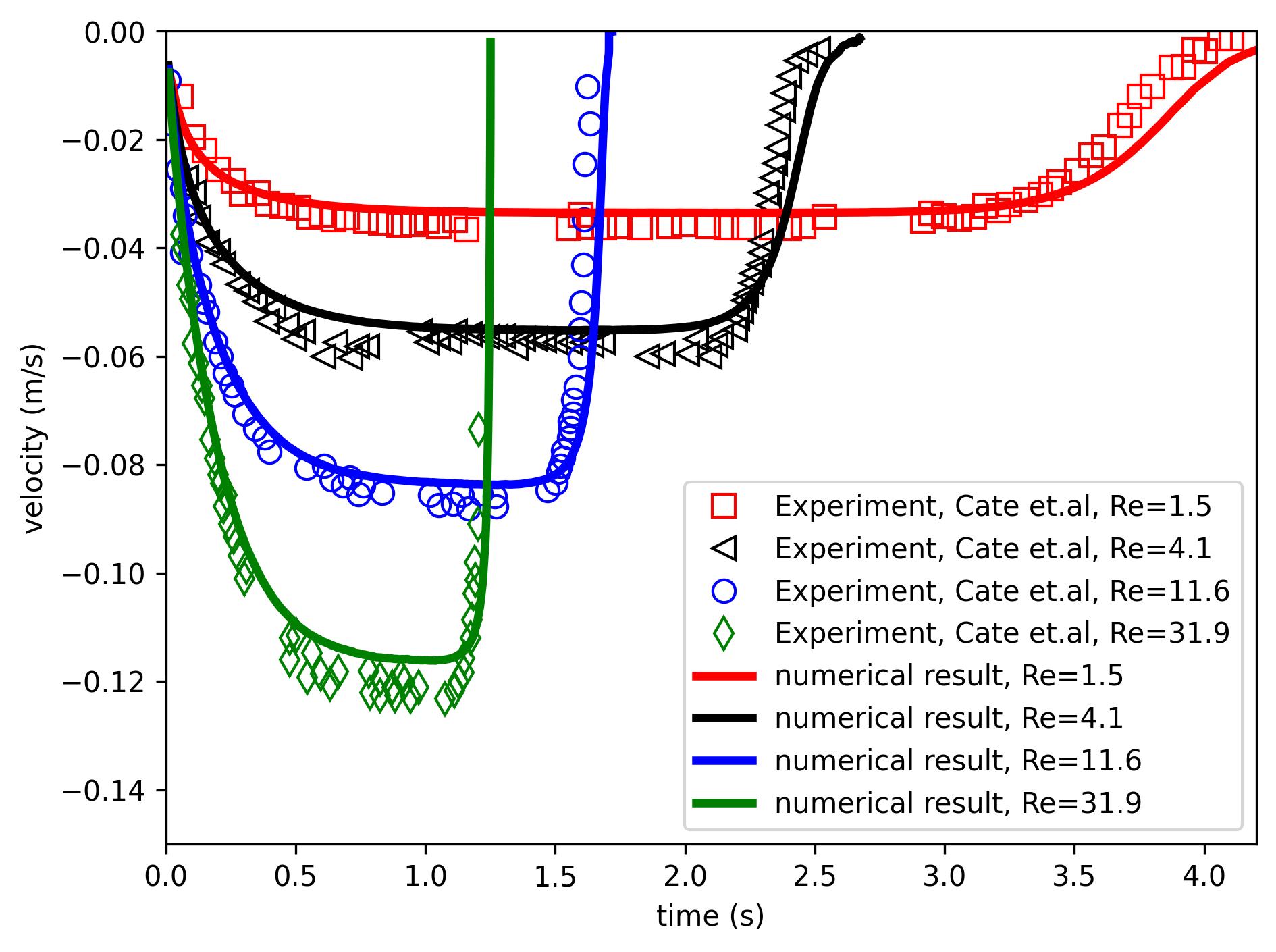}
    \caption{Falling velocity vs time}
  \end{subfigure}
  \caption{(a) Time evolution of vertical position presented as height-diameter ratio. (b) Time evolution of falling velocity. Experimental data from \cite{ten2002particle} is shown by symbols in the figure.}
  \label{fig:freefall_sketch}
\end{figure}

{Figure \ref{fig:freefall_different_time} shows the simulated velocity field of case 1 at four time instants. The color encodes the magnitude of the velocity for both the sphere and the surrounding fluid, and the arrows near the sphere indicate the flow direction. Figures \ref{fig:freefall_sketch}(a) and (b) show the time evolution of the sphere height-diameter ratio and the falling velocity for the four simulated cases, with the experimental data of \cite{ten2002particle} overlaid for comparison.}

{The simulated height-diameter ratio and the falling velocity match the experimental data of \cite{ten2002particle} across the four viscosity cases, with deviations remaining within the experimental uncertainty. The match holds across the Reynolds number range $Re \in [1.5, 31.9]$, confirming that the three-dimensional implementation reproduces the inertial-viscous balance during the gravitational acceleration phase and the asymptotic terminal-velocity regime.}


\subsection{Turek--Hron FSI3 benchmark}
\begin{figure}[h]
  \centering
  \includegraphics[width=0.95\textwidth]{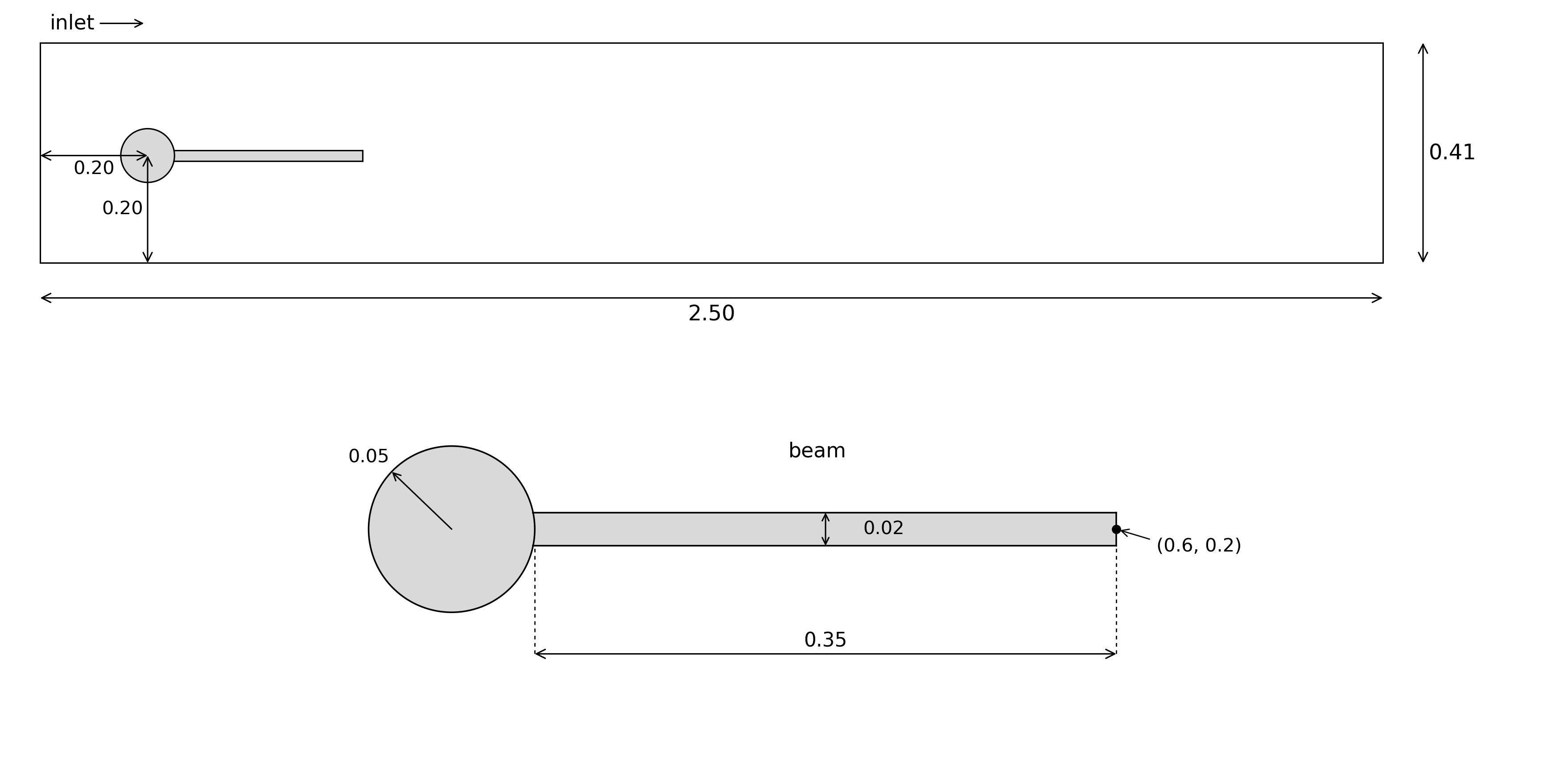}
  \caption{Geometric setup of the Turek--Hron FSI3 benchmark with the beam-tip tracking point $(0.6,0.2)$ used for response comparison.}
  \label{fig:turek_fsi3_setup}
\end{figure}

The Turek--Hron FSI3 benchmark couples a laminar channel flow past a fixed circular cylinder with an elastic beam attached downstream. The matched fluid-solid density ($\rho_s/\rho_f = 1$) makes FSI3 a stringent stability probe: classical staggered partitioned schemes diverge under the strong added-mass effect, and the test directly assesses the ALE solver's behavior in the regime where added-mass coupling, large structural deformation, and moving-mesh stability act simultaneously. The benchmark geometry is shown in Figure \ref{fig:turek_fsi3_setup}, and the parameters used in this work are summarized in Table \ref{table:turek_fsi3_setup}.

\begin{table}[h]
  \centering
  \caption{Flow and material parameters of the Turek--Hron FSI3 benchmark.}
  \label{table:turek_fsi3_setup}
  \begin{tabular}{l c}
    \toprule
    Parameter & Value \\
    \midrule
    Fluid density $\rho_f$ & $1000$ \\
    Fluid kinematic viscosity $\nu_f$ & $1.0\times10^{-3}$ \\
    Mean inlet velocity & $2$ \\
    Reynolds number & $200$ \\
    Solid density $\rho_s$ & $1000$ \\
    Solid Poisson ratio $\nu_s$ & $0.4$ \\
    Solid elastic parameter $\mu_s$ & $2.0\times10^{6}$ \\
    \bottomrule
  \end{tabular}
\end{table}

\begin{figure}[h]
  \centering
  \includegraphics[width=0.95\textwidth]{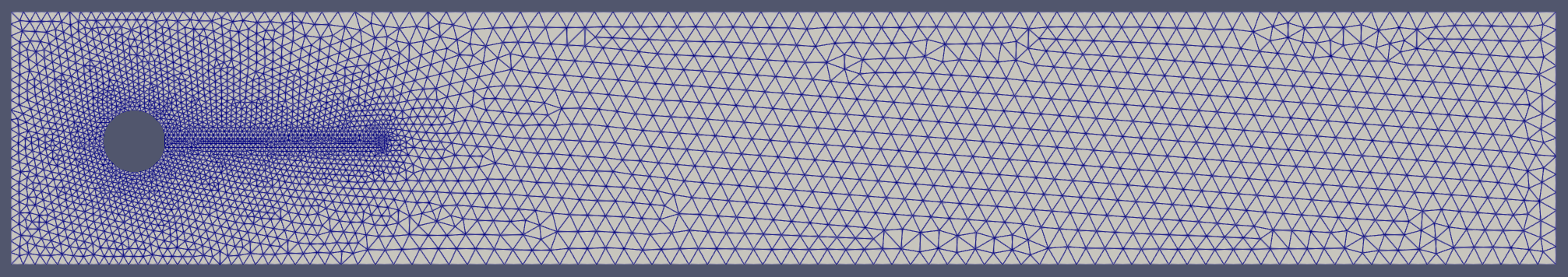}
  \caption{Computational mesh used for the Turek--Hron FSI3 benchmark in this work.}
  \label{fig:turek_fsi3_mesh}
\end{figure}

The inlet velocity is prescribed as a parabolic channel profile with a time-ramped mean velocity. Specifically, the inflow condition is given by
\[
\boldsymbol{u}_{\mathrm{in}}(y,t)
=
\left(
\frac{6\bar{U}(t)}{H^2}y(H-y),0
\right),
\qquad
\bar{U}(t)=
\begin{cases}
t, & 0\le t < 2,\\
2, & t\ge 2.
\end{cases}
\]
Thus, the mean inlet velocity increases linearly from zero to its maximum value of $2$ at $t=2$, and remains constant afterwards.

In our implementation, the downstream end of the beam is rounded to reduce the local mesh distortion and numerical oscillations caused by the sharp trailing corner. This small geometric regularization preserves the global FSI response while improving the robustness of the ALE mesh motion near the beam tip. For quantitative comparison, we also track the vertical and horizontal responses of the center point at the downstream end of the beam.

\begin{figure}[h]
  \centering
  \includegraphics[width=0.95\textwidth]{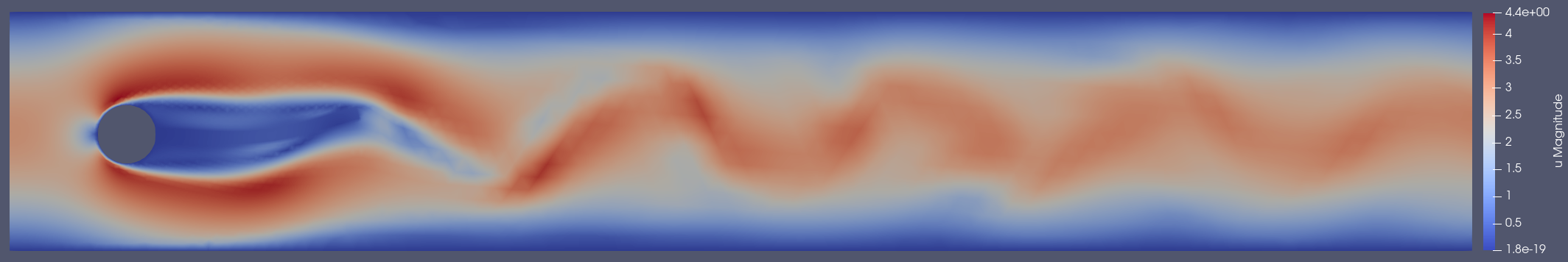}
  \caption{Velocity magnitude distribution for the Turek--Hron FSI3 benchmark at the instant when the beam reaches its maximum displacement in the $y$-direction.}
  \label{fig:turek_fsi3_velocity}
\end{figure}

\begin{figure}[!htbp]
  \centering
  \begin{subfigure}{0.8\textwidth}
    \centering
    \includegraphics[width=\textwidth]{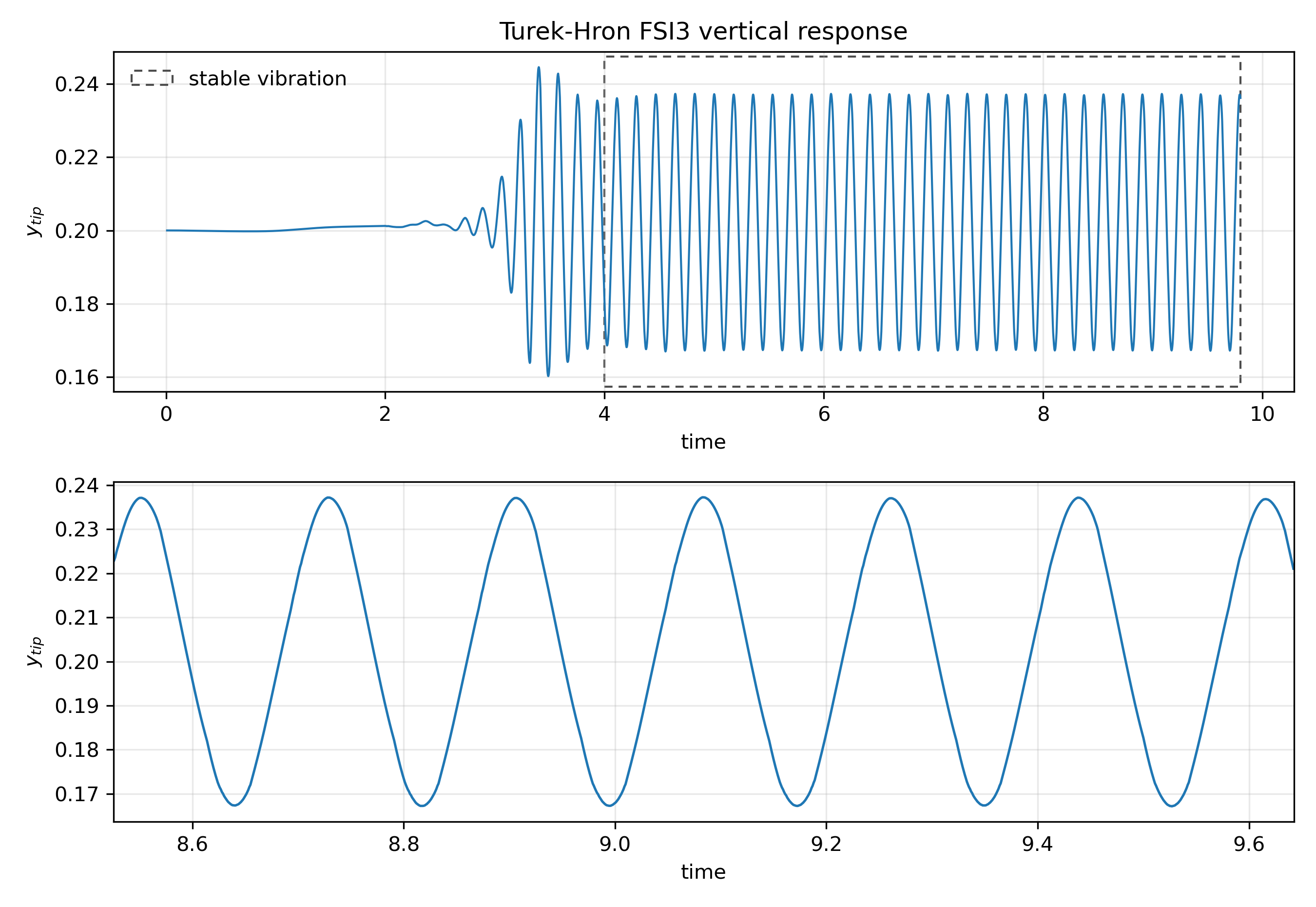}
    \caption{Vertical response of the beam-tip tracking point.}
    \label{fig:turek_fsi3_y_response}
  \end{subfigure}

  \begin{subfigure}{0.8\textwidth}
    \centering
    \includegraphics[width=\textwidth]{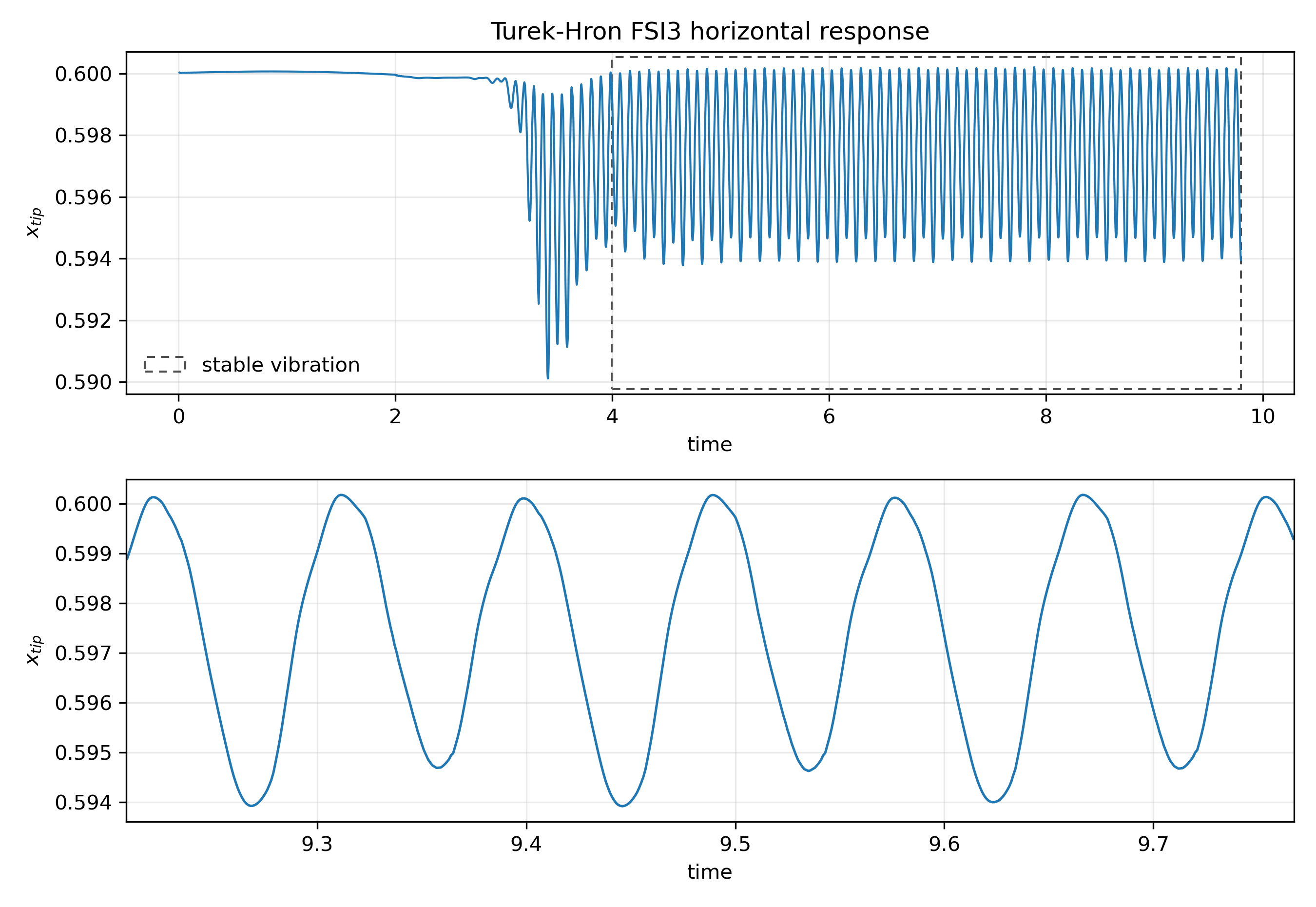}
    \caption{Horizontal response of the beam-tip tracking point.}
    \label{fig:turek_fsi3_x_response}
  \end{subfigure}
  \caption{Beam-tip response histories for the Turek--Hron FSI3 benchmark.}
  \label{fig:turek_fsi3_tip_response}
\end{figure}

The response histories in Figure \ref{fig:turek_fsi3_tip_response} show a stable periodic oscillation after the initial transient. The computed $y$-direction vibration frequency is $5.6333$ with amplitude $0.0348$, while the $x$-direction frequency is $11.2656$, consistent with the expected frequency-doubling between streamwise and transverse responses driven by the symmetric vortex shedding. The reference values reported by Turek and Hron are $5.3$ and $0.03438$, giving relative discrepancies below $3\%$ in both amplitude and frequency. The matched-density configuration imposes the strongest added-mass coupling in the benchmark series, and the stable convergence of the quasi-monolithic solver under this regime confirms that the implicit fluid-solid block is sufficient to bypass the artificial added-mass instability that destabilizes staggered partitioned schemes. Furthermore, the IMEX-PRK scheme preserves the rapid elastic transients without the dissipative damping characteristic of backward Euler, which is critical for accurately capturing the limit-cycle amplitude.

\subsection{Multi-scale simulation in long spiral channel}
Recently, as the number of reports on microfluidic devices is increasing fast, many studies have focused on the focusing and separation of particles in microchannels with spiral structures \cite{chen2025microfluidic, feng2022viscoelastic, kumar2021high}. In some researches, obstacles are added inside the spiral channel in order to improve their performance \cite{zhao2022flow, shen2024spiral}. In this part, we will use the local updating algorithm to present the simulation results for these cases.

\begin{figure}[h]
  \centering
  \includegraphics[width=0.6\textwidth]{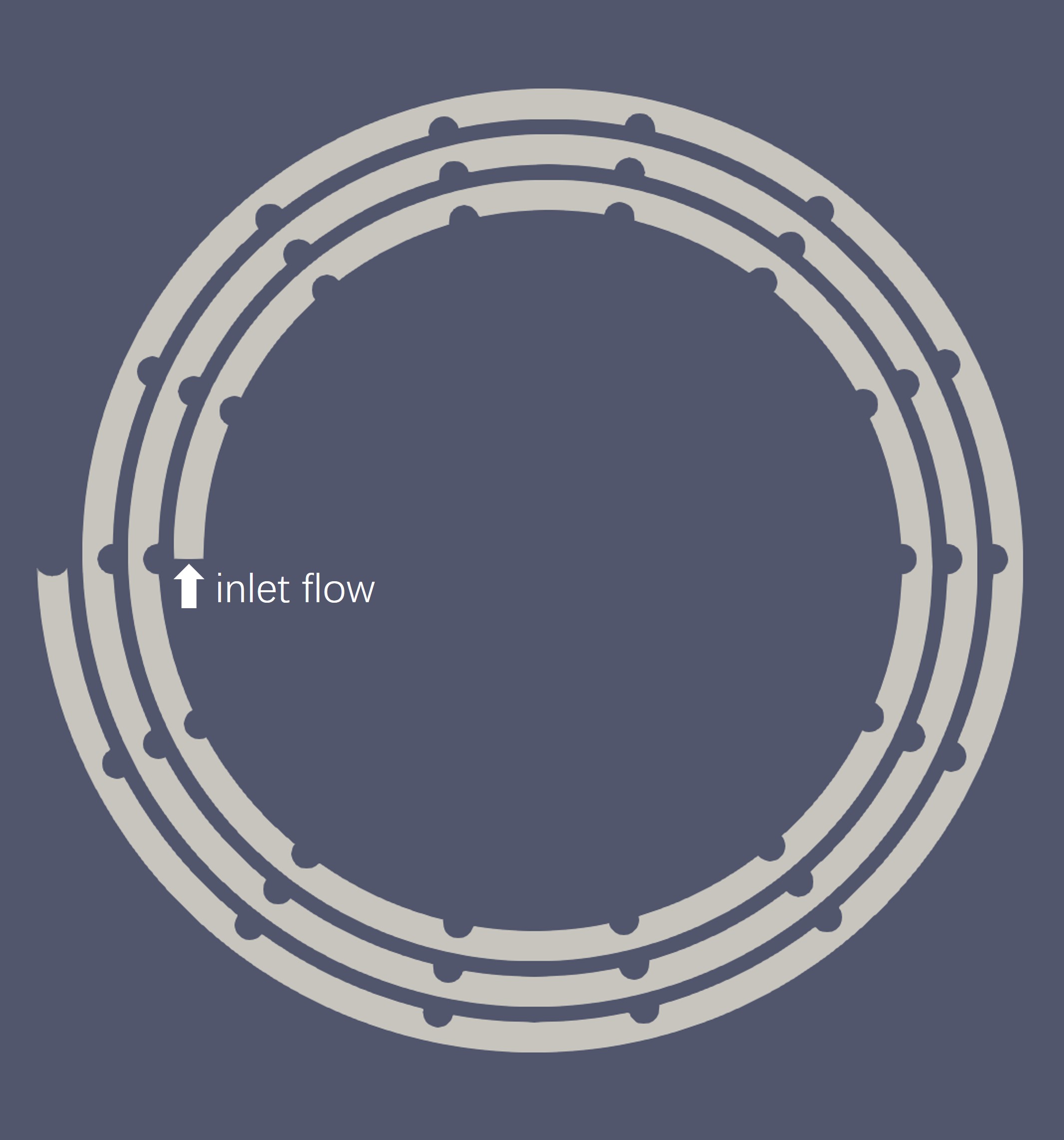}
  \caption{The background domain of the spiral channel.}
  \label{fig:spiral}
\end{figure}

As shown in Figure \ref{fig:spiral} a spiral channel with obstacles are set. According to the reference above, the radius of the innermost circular arc is $10000\mu m$, which is typical scale in recent researches. The radius of the arc gradually increases, with each loop increasing by $1350\mu m$. The channel is $900\mu m$ wide and $100\mu m$ high. Each loop contains 14 obstacles with radius of $450\mu m$. The flow enters through the inner port and exits through the outer port. At the entrance velocity Dirichlet boundary condition is applied to ensure the flux of the channel is adjustable.

\begin{figure}[!htbp]
  \centering
  \begin{subfigure}{0.48\textwidth}
    \centering
    \includegraphics[height=4.5cm]{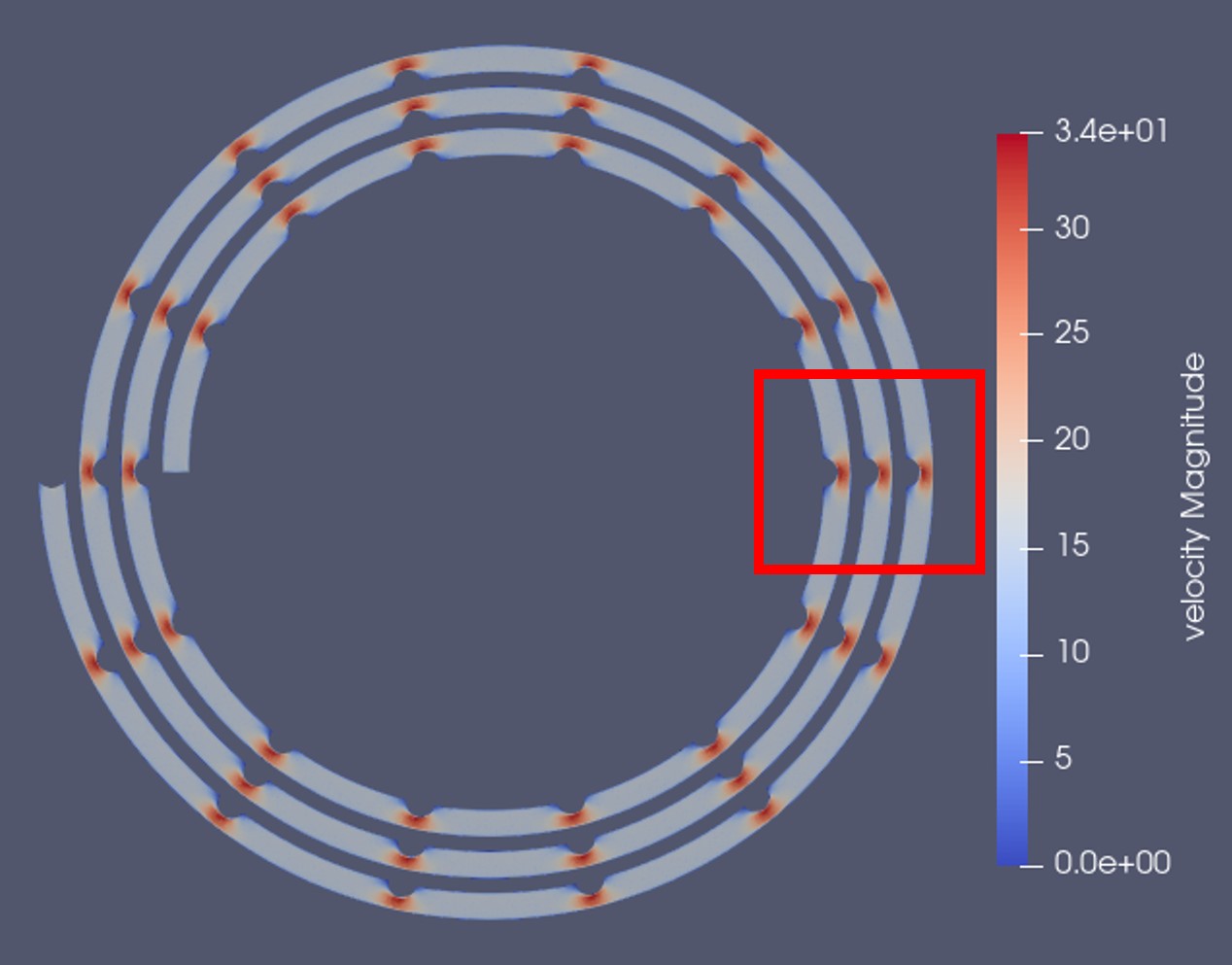}
    \caption{flux=1000$\mu$L/min}
  \end{subfigure}\hfill
  \begin{subfigure}{0.48\textwidth}
    \centering
    \includegraphics[height=4.5cm]{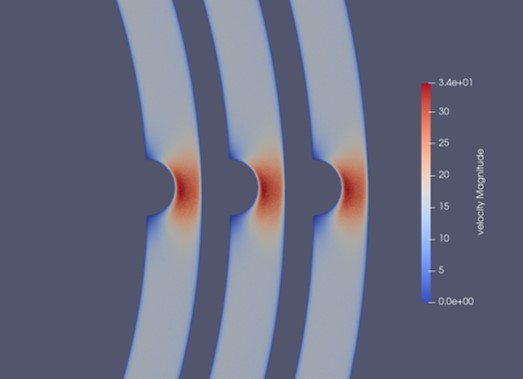}
    \caption{zoomed view of flux=1000$\mu$L/min}
  \end{subfigure}

  \begin{subfigure}{0.48\textwidth}
    \centering
    \includegraphics[height=4.5cm]{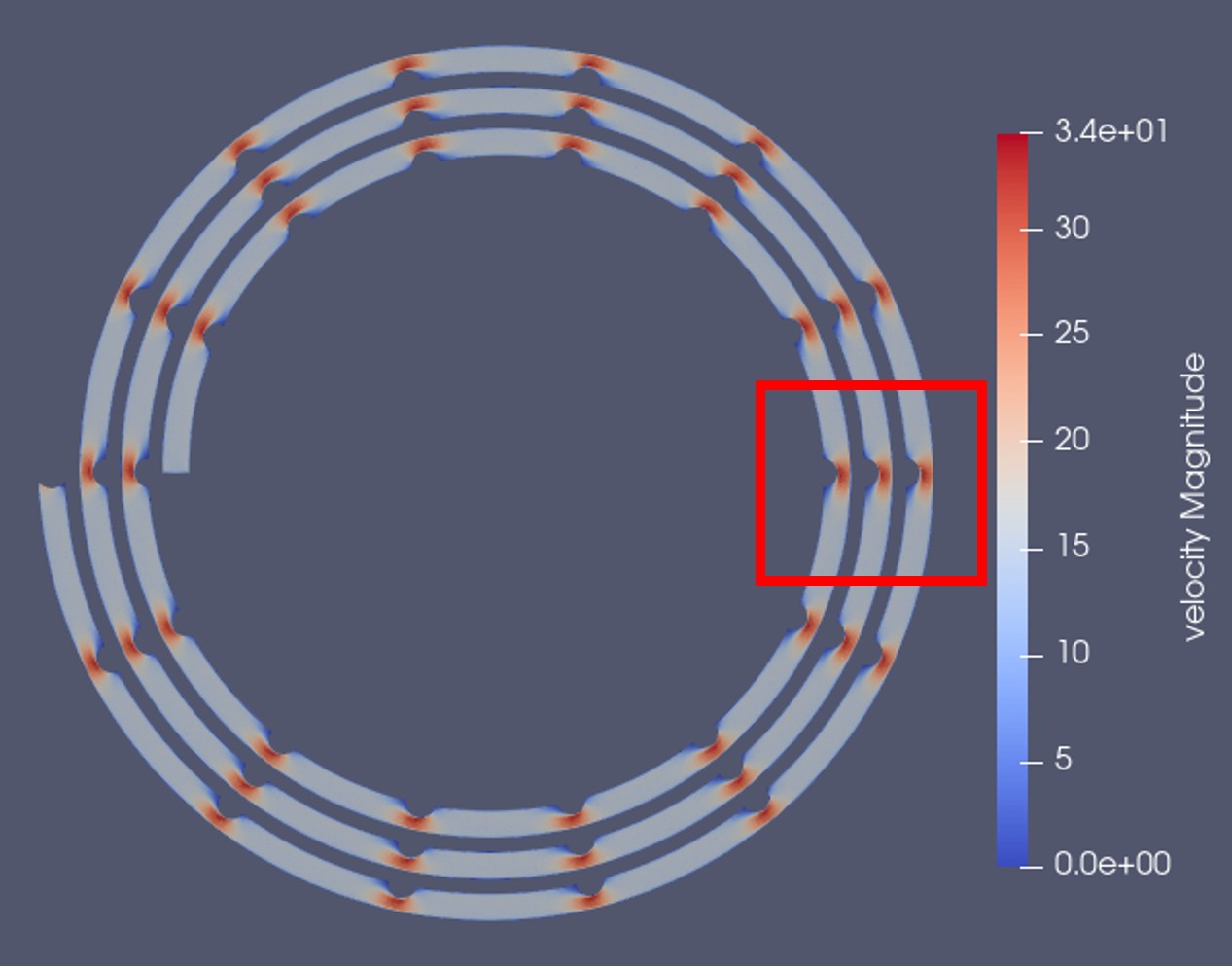}
    \caption{flux=1500$\mu$L/min}
  \end{subfigure}\hfill
  \begin{subfigure}{0.48\textwidth}
    \centering
    \includegraphics[height=4.5cm]{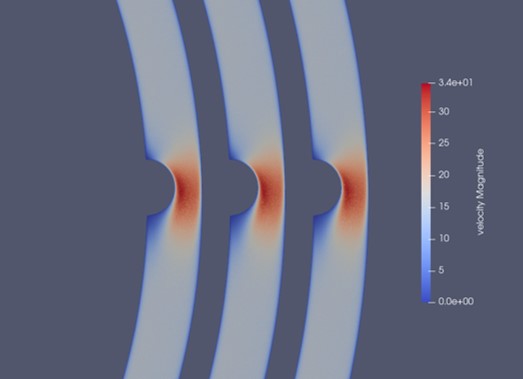}
    \caption{zoomed view of flux=1500$\mu$L/min}
  \end{subfigure}

  \begin{subfigure}{0.48\textwidth}
    \centering
    \includegraphics[height=4.5cm]{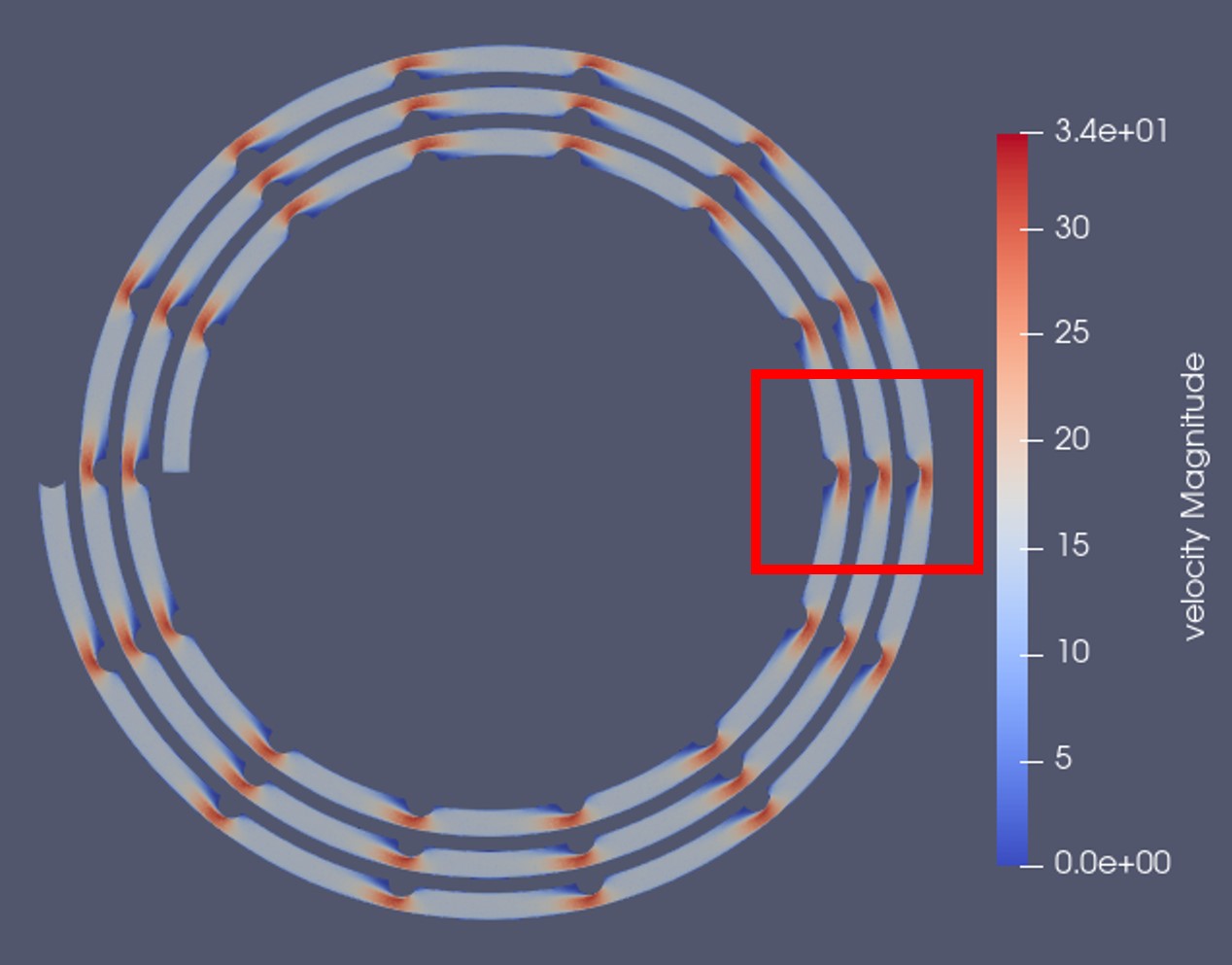}
    \caption{flux=3000$\mu$L/min}
  \end{subfigure}\hfill
  \begin{subfigure}{0.48\textwidth}
    \centering
    \includegraphics[height=4.5cm]{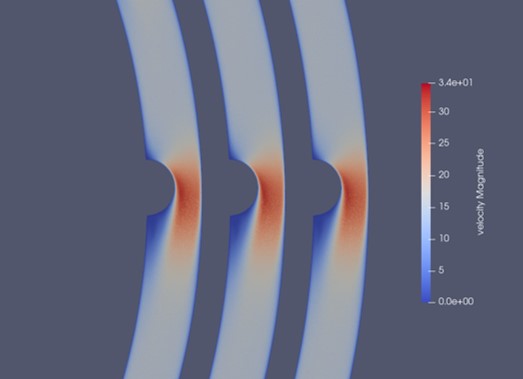}
    \caption{zoomed view of flux=3000$\mu$L/min}
  \end{subfigure}

  \caption{The magnitude of the velocity field on the cross-section of the flow channel with different flux (left) along with the zoomed views (right) of the local part in the red frame.}
  \label{fig:spiral_velocity}
\end{figure}

Figure \ref{fig:spiral_velocity} plotted the magnitude of velocity field at height $50 \mu m$ with flux $1mL/min$, $1.5mL/min$ and $3mL/min$. Particles with radius $r=4$ are released in the channel with 0 velocity initially at the inner loop, distributed uniformly along y-axis with $z=-25\mu m$ and $x=0\mu m$.

\begin{figure}[htbp]
  \centering
  \begin{subfigure}{0.8\textwidth}
    \includegraphics[width=\linewidth]{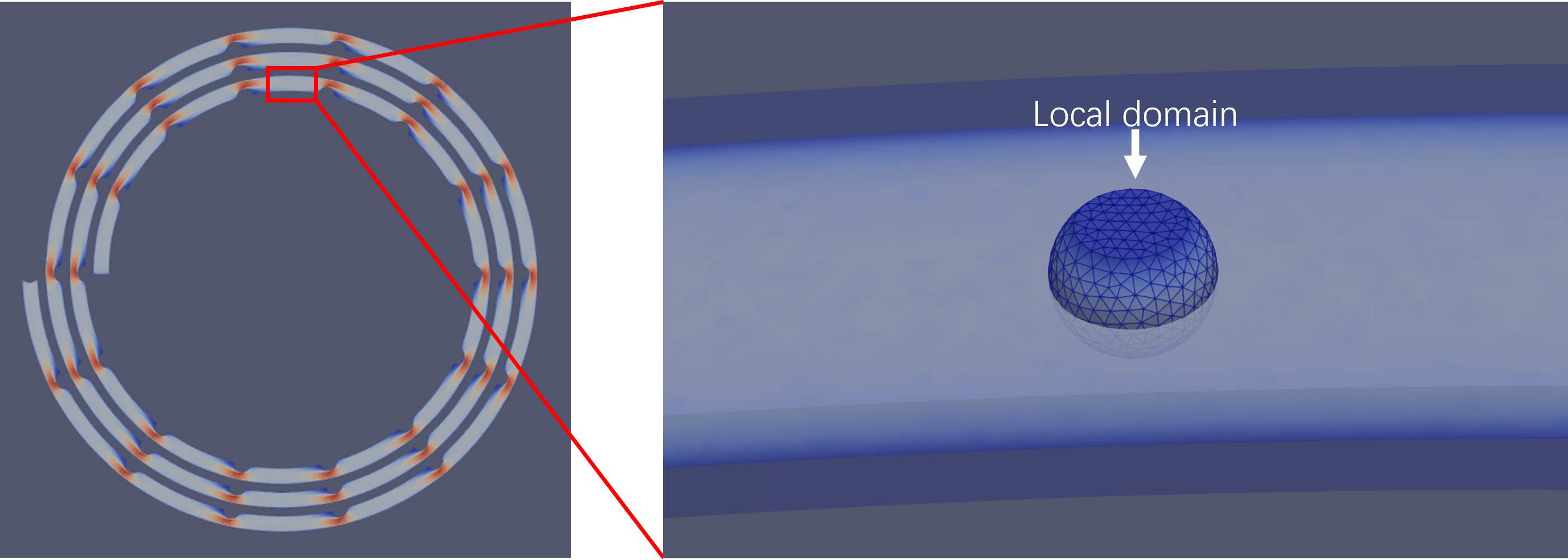}
    \caption{Local domain extraction}
  \end{subfigure}

  \vspace{1.5em} 

  \begin{subfigure}{0.8\textwidth}
    \includegraphics[width=\linewidth]{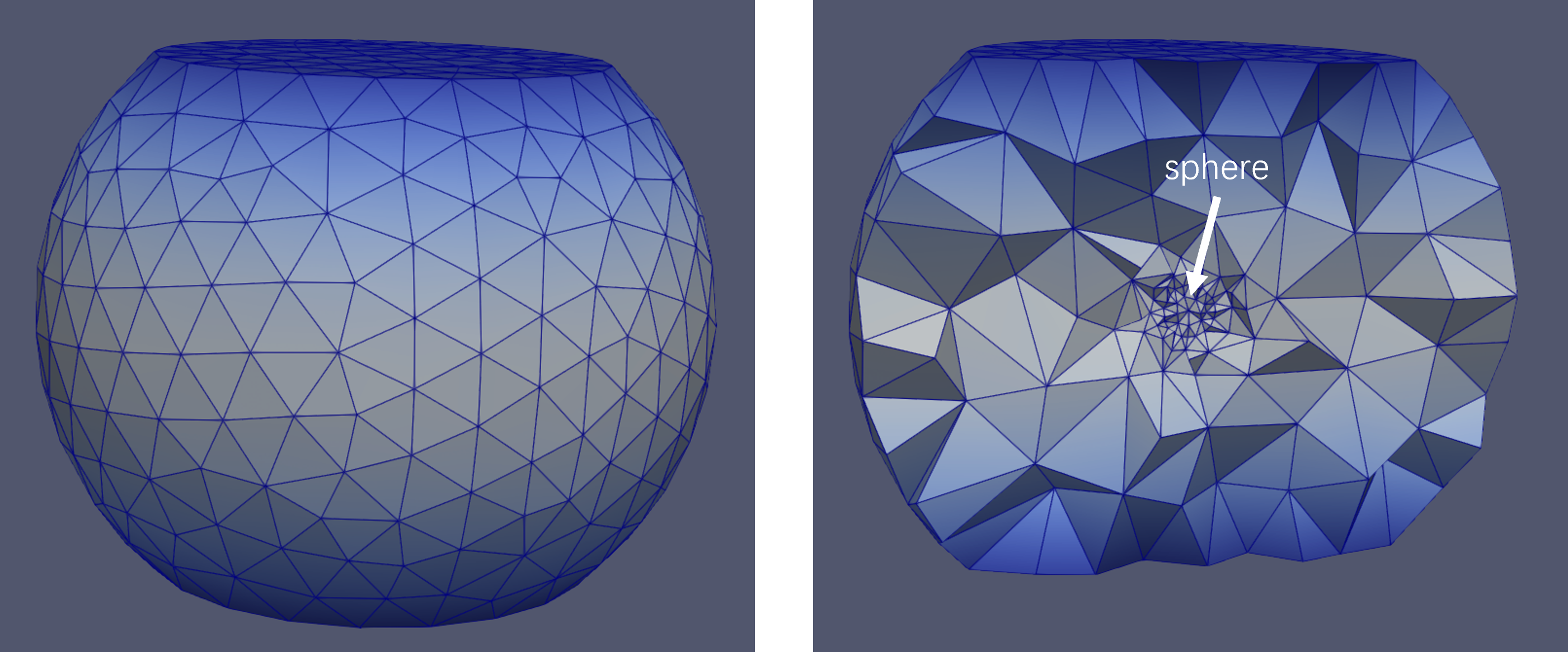}
    \caption{Local domain and cross section together with velocity field}
  \end{subfigure}

  \caption{(a) Local domain extraction. As can been seen, there is a significant scale disparity between the flow channel and the local domain. (b) The mesh of local domain together with the cross section. The sphere is located in the center of the domain and the mesh size nearby is thinner. The calculated velocity field is also plotted on the meshes.}
  \label{fig:spiral_local_domain}
\end{figure}

The FSI calculation is done by the local update algorithm in order to improve the efficiency of the simulation. The local domain extracted from the entire flow channel is shown in Figure \ref{fig:spiral_local_domain} and presented together with the calculated velocity field.

\begin{figure}
  \centering
  \begin{subfigure}{\textwidth}
    \centering
    \includegraphics[width=0.66\textwidth]{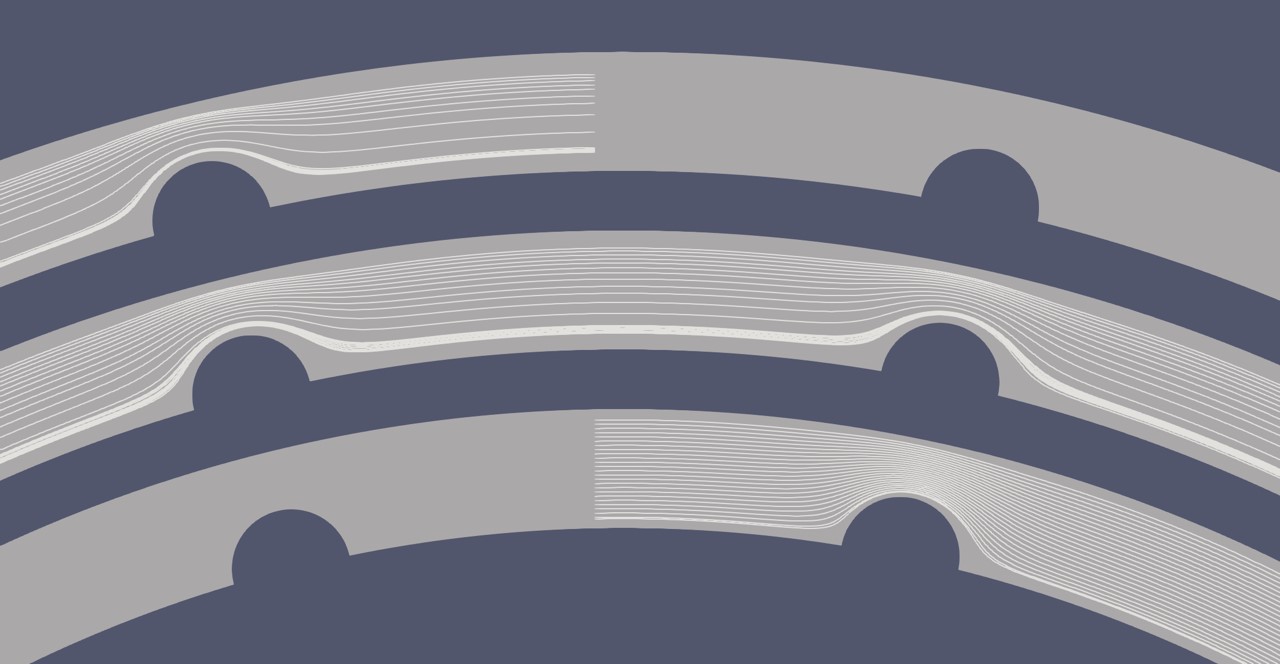}
    \caption{The trajectories of particles released at with flux=1000$\mu$L/min}
  \end{subfigure}
  \begin{subfigure}{\textwidth}
    \centering
    \includegraphics[width=0.66\textwidth]{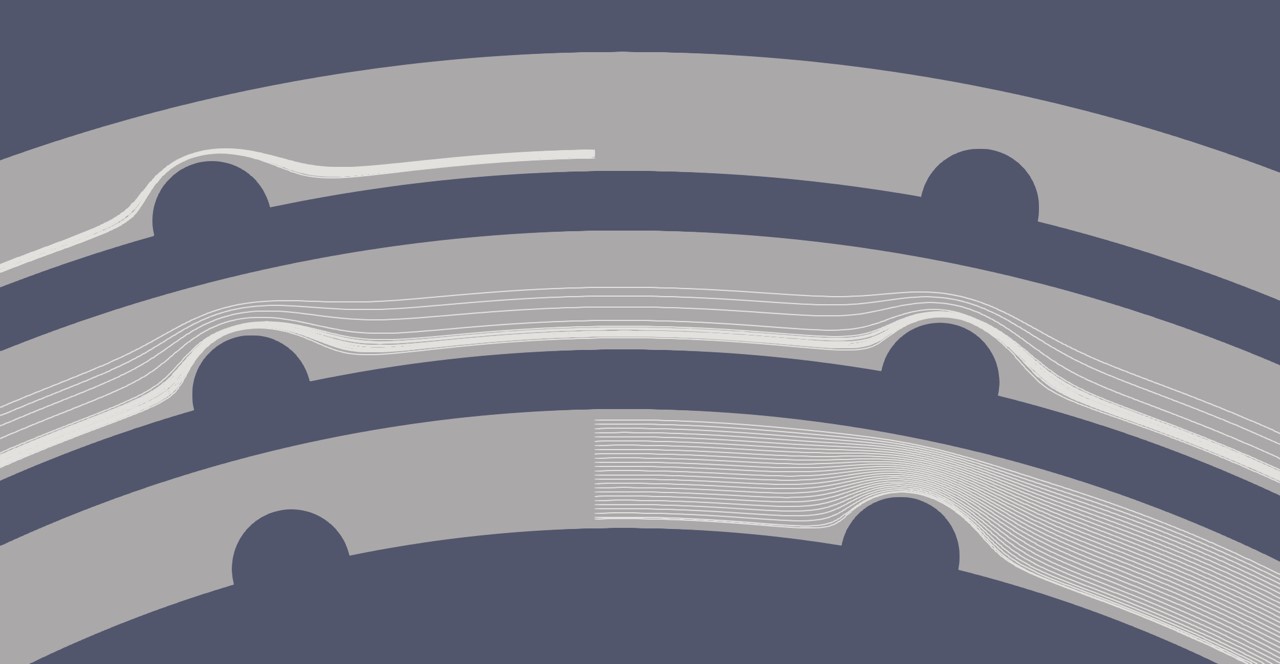}
    \caption{The trajectories of particles released at with flux=1500$\mu$L/min}
  \end{subfigure}
  \begin{subfigure}{\textwidth}
    \centering
    \includegraphics[width=0.66\textwidth]{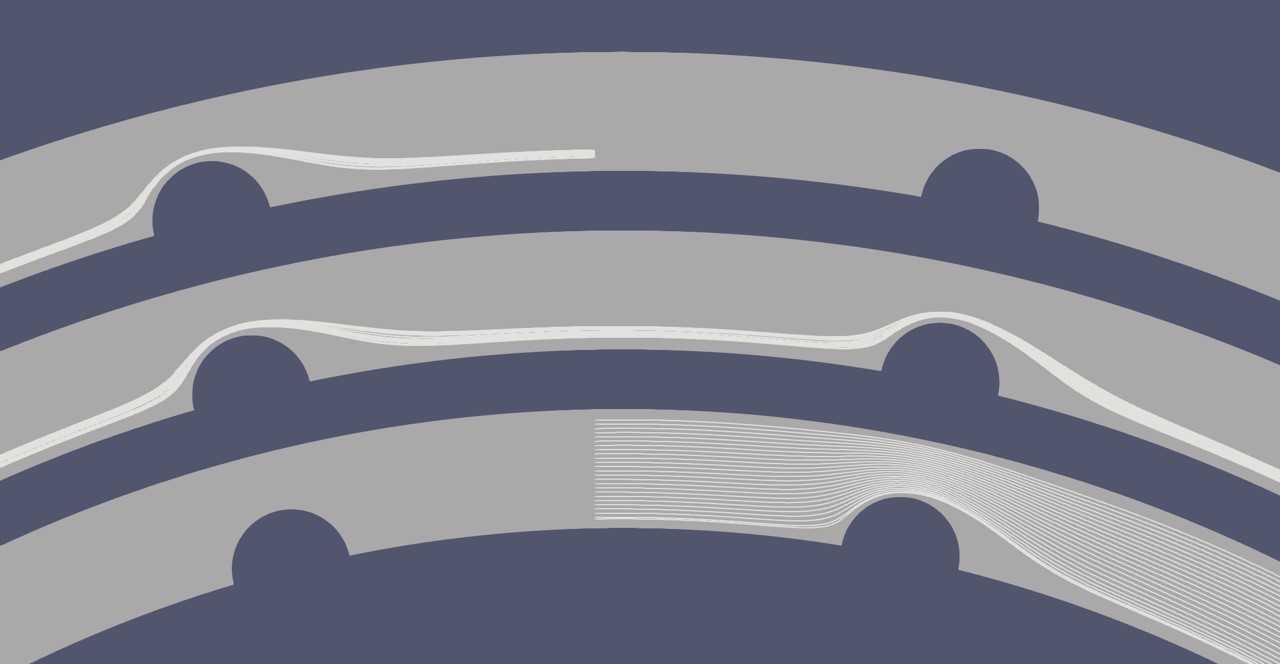}
    \caption{The trajectories of particles released at with flux=3000$\mu$L/min}
  \end{subfigure}

  \caption{25 particles are released at the inner loop of the channel at with $z=-25\mu m$. Two loops are traveled totally, as can be seen from the trajectories, the focusing phenomenon is increasingly clear as the flux in the channel increase.}
  \label{fig:spiral_trajectory}
\end{figure}

Two loops are simulated and the trajectories of the particles are shown in Figure \ref{fig:spiral_trajectory}. As the flow rate increases, particle focusing along the spiral channel is progressively enhanced. The trend aligns with established theory: an elevated flow rate, characterized by an increased Dean number $De$, strengthens the secondary flow induced by the curvature. The centrifugal pressure gradient across the curved cross-section intensifies the counter-rotating Dean vortices, which exert a more potent Dean drag force $F_D$ on the suspended particles. The balance between the inertial lift force $F_L$ and the intensified Dean drag is shifted toward the equilibrium position, accelerating the convergence of the particle ensemble to a narrow focusing band. The presence of obstacles inside the spiral channel introduces additional localized hydrodynamic perturbations: the geometric constraints locally modulate the secondary-flow pattern and generate high-shear zones that refine the focusing width. These observations confirm that the localized updating strategy preserves the relevant hydrodynamic mechanisms across both the local FSI scale and the channel-scale advection.



\section{Conclusion}
 {This work develops the qMLH-ALE framework: a quasi-monolithic localized high-order finite element method for multi-scale fluid-structure interaction, combining a localized updating algorithm for multi-scale simulation with high-order temporal and spatial discretizations on the localized sub-problem.}

 {
  Three contributions support the framework's applicability. The sharp-interface model on second-order body-fitted grids attains third-order geometric approximation of curved boundaries and removes the interfacial smearing characteristic of fixed-grid methods. The use of isoparametric elements maintains the optimal spatial convergence rate, which directly improves the accuracy of boundary stress and wall-shear computations relevant to aerospace and naval applications. The second-order IMEX-PRK time integration on the moving mesh achieves second-order temporal accuracy while avoiding the dissipative damping of backward Euler, making the algorithm suitable for tracking strongly transient interfaces such as red blood cells in complex flows.

  The Turek-Hron FSI3 benchmark, performed under matched fluid-solid density, agrees with the reference data within $3\%$ in beam-tip amplitude and frequency, providing self-proof of the solver's stability under the strong added-mass effect that destabilizes conventional partitioned schemes. Simulations of particle focusing in spiral microchannels reproduce the experimentally observed flux-dependent focusing behavior, confirming that the localized updating strategy preserves the relevant hydrodynamic mechanisms while reducing the cost of long-range multi-scale FSI to a tractable level. The framework thus provides a computational basis for the design of microfluidic devices and other multi-scale engineering systems where geometric fidelity and long-range tracking must be retained simultaneously.
 }

\section*{Declaration of generative AI and AI-assisted technologies in the manuscript preparation process}
During the preparation of this work the authors used Gemini in order to improve the clarity of the manuscript. After using this tool/service, the authors reviewed and edited the content as needed and take full responsibility for the content of the published article.

\section*{Declaration of interest}
The authors declare that they have no known competing financial interests or personal relationships that could have appeared to influence the work reported in this paper.

\FloatBarrier
\newpage
\printbibliography

\end{document}